\documentclass[12pt]{article}
\usepackage{amsfonts}
\usepackage{epsfig}
\title{Introducing The Plaid Model}
\author{Richard Evan Schwartz \thanks{\hskip 5 pt Supported by 
N.S.F. Research Grant DMS-1204471}}

\newtheorem{theorem}{Theorem}[section]

\newtheorem{lemma}[theorem]{Lemma}

\newtheorem{corollary}[theorem]{Corollary}
\newtheorem{conjecture}[theorem]{Conjecture}

\def\startproof{{\bf {\medskip}{\noindent}Proof: }}

\def\endproof{$\spadesuit$  \newline}

\def\Q{\mbox{\boldmath{$Q$}}}% 
\def\R{\mbox{\boldmath{$R$}}}% 
\def\Z{\mbox{\boldmath{$Z$}}}% 

\begin{document}
\maketitle

\begin{abstract}
We introduce and prove some basic
results about a combinatorial model which
produces embedded polygons in
the plane.  The model is highly
structured and relates to a
number of topics in dynamics and geometry.
\end{abstract}

\section{Introduction}

The purpose of this paper is to introduce a combinatorial
model which produces embedded lattice polygons in the
plane.  The model is related to a variety of things:
outer billiards on kites [{\bf S1\/}], [{\bf S2\/}],
circle rotations, polyhedron exchange
transformations, corner percolation, 
P. Hooper's Truchet
tile system [{\bf H\/}], and 
DeBruijn's theory of multigrids [{\bf DeB\/}].
I call the model the {\it plaid model\/} because of
the grids of parallel lines it involves.

The plaid model depends on a rational parameter
$p/q \in (0,1)$, and we always \footnote{There is a
similar kind of theory when $pq$ is odd, but it
is sufficiently different that I will not discuss it here.
Also, I haven't really worked it out in detail.}
take $pq$ to be
even.  Figure 1.1 shows some of the polygons
produced by the model for the parameter $4/17$
and Figure 1.2 shows some of the polygons
produced by the model for the parameter
$17/72$.  These parameters
are successive terms in the continued
fraction expansion of $\sqrt 5-2$, a 
parameter which figures heavily in [{\bf S2\/}].
If you look closely at these pictures
you notice that the second one contains lots of copies
of the polygons in the first one, and also the
large polygons in the two pictures
match up very well when they are rescaled to
have the same size and superimposed.

\newpage
\begin{center}
\resizebox{!}{2.2in}{\includegraphics{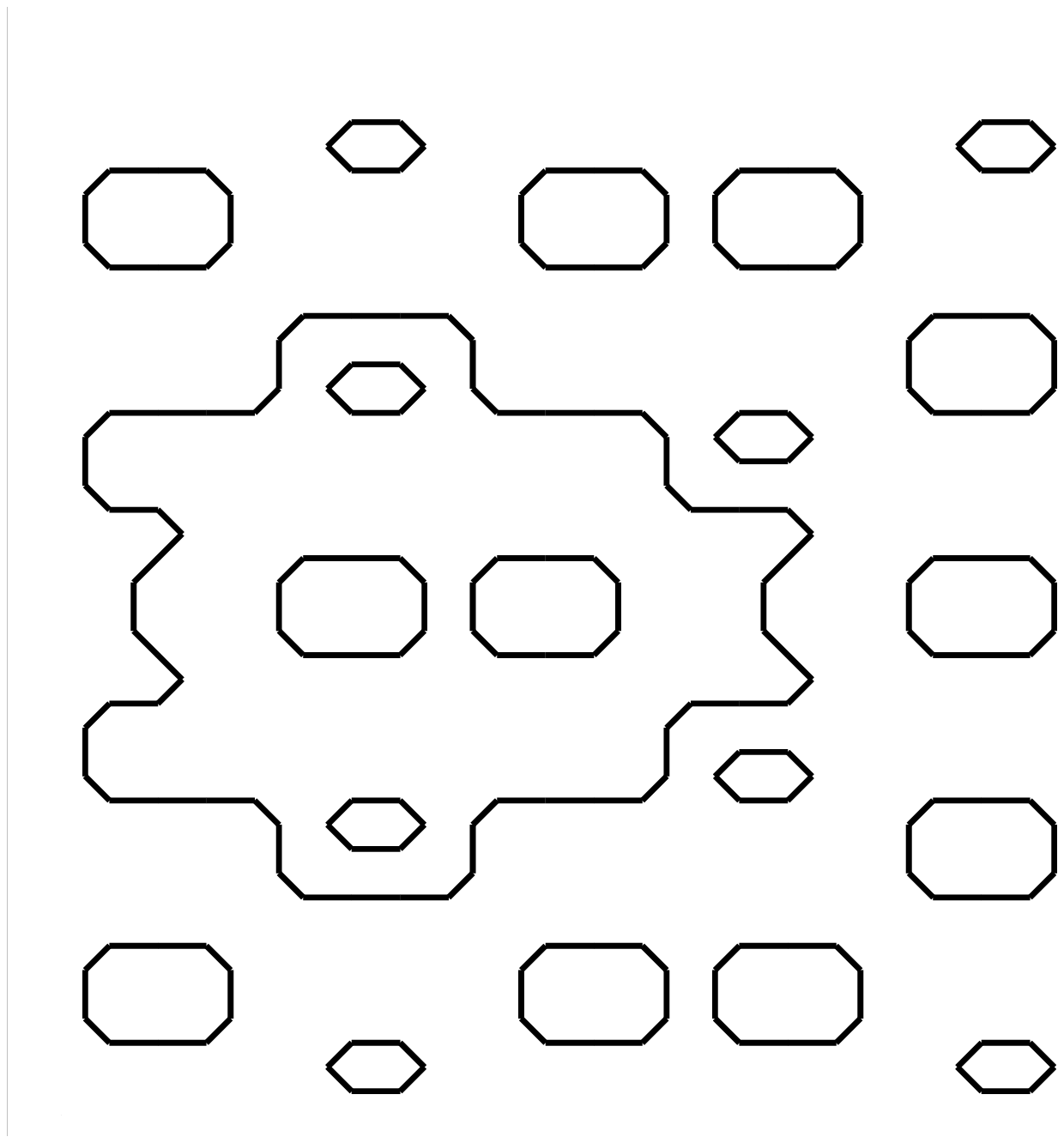}}
\newline
{\bf Figure 1.1:\/} Some of the polygons
for the parameter $4/17$.
\end{center}

\begin{center}
\resizebox{!}{4.8in}{\includegraphics{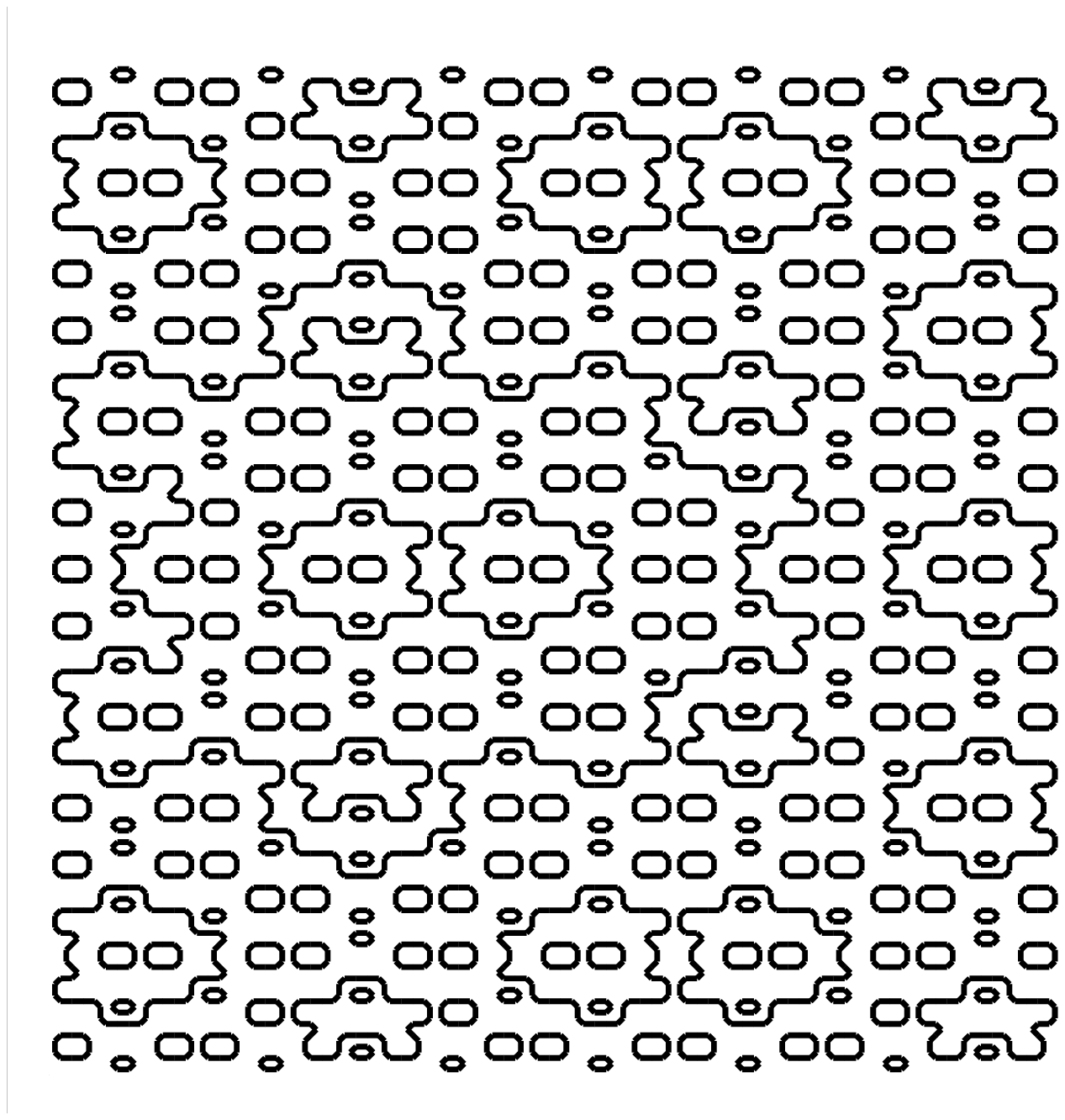}}
\newline
{\bf Figure 1.2:\/} Some of the polygons
for the parameter $17/72$.
\end{center}

These kinds of phenomena are ubiquitous in
the model, and hint at its depth.
My main aim in this paper is to introduce
the plaid model and prove some basic
things about it.

I discovered the plaid model in my effort to
understand the coarse self-similarity
one sees in the {\it arithmetic graphs\/}
associated to outer billiards on kites.
There seems to be so much to say about
the plaid model that I don't even get
to the outer billiards connection
in this paper, though I will include a
brief account in \S 8.  
The plaid model is a kind of elaboration of
the Hexagrid Theorem from [{\bf S1\/}] but
the account here is self-contained
and has nothing apparently to do with
outer billiards on kites.
I imagine that the plaid model is something that
one could study for its own sake.

Here is an overview of the paper.
In \S 2, I will describe the plaid model
in terms of the intersection points of a system
of grids in the plane, and I will establish
some of its basic properties.
From the description, it is not at all
clear why the model produces polygonal
paths, or indeed anything interesting.
However, the Theorem \ref{fun1} (the
Fundamental Theorem)
says that the plaid model really
does make sense as a generator of
polygonal paths. 

In \S 3 I will give a second description
of the plaid model in terms of
polyhedron exchange transformations.
Here is the main result of the paper. 

\begin{theorem}[PET Equivalence]
\label{tile}
Let $p/q$ be an even rational parameter.
Let $P=2p/(p+q)$.  The plaid polygons
associated to $p/q$
describe the vector dynamics of a set
of distinguished orbits in a
$3$ dimensional polyhedron exchange 
transformation $\widehat X_P$.
Moreover, there is a $4$-dimensional
fibered convex integral affine polytope
exchange transformation whose
fiber over $P$ is
$\widehat X_P$.
\end{theorem}

All this terminology will be defined
in \S 3.
By {\it vector dynamics\/}, I mean that
we associate one of the unit
vectors in $\Z^2$, or the $0$-vector,
to each region 
of a polyhedron exchange transformation.
Following an orbit, one obtains a
sequence of vectors, which one ``accumulates''
to produce a polygonal path.  This will
make sense even at irrational parameters
but when $P=2p/(p+q)$ the vector
dynamics of the distinguished orbits
produces the plaid polygons with
respect to the paramater $p/q$.

For each parameter, our PET has a
flat $3$ torus as its domain.
This flat $3$ torus has a natural
subdivision into cubes of sidelength
$2/(p+q)$, and the dynamics preserves the
set of centers of these cubes.  These
are the distinguished orbits referred
to in the PET Equivalence Theorem.

  As we
will see, the PET Equivalence Theorem immediately
implies \footnote{Technically, we will deduce the
Fundamental Theorem from Theorem \ref{iso}, the
Isomorphism Theorem, which is a precursor
to the PET Equivalence Theorem.}
 the Fundamental Theorem. 
 Moreover,
the PET Equivalence Theorem allows one to extend
the plaid model to irrational parameters in
$(0,1)$.

In \S 4 I will prove several results concerning
the distribution of large orbits in the
rational case.  Theorems \ref{bigg} and
\ref{bigg2} are the main results.
Theorem \ref{bigg2} goes
beyond anything I could say about
outer billiards in [{\bf S1\/}],
and involves a rescaled limit
of the plaid model. 

In \S 5-8 I will prove the PET Equivalence Theorem.
The proof exploits the
``physical nature'' of the plaid
model, in which certain of the intersection
points are grouped into ``particles'' which
(when suitably interpreted) move around and
obey billiard-like laws.

In \S 9 I will explain the connection
between the plaid model and outer billiards
on kites.   The main conjectural result,
Conjecture \ref{qi} (The Quasi-Isomorphism
Conjecture)
says that the closed loops produced
by the plaid model are the same,
up to an error of at most $2$ units,
as certain affine images of the arithmetic
graphs associated
to outer billiards on kites.  

\begin{center}
  \resizebox{!}{2.8in}{\includegraphics{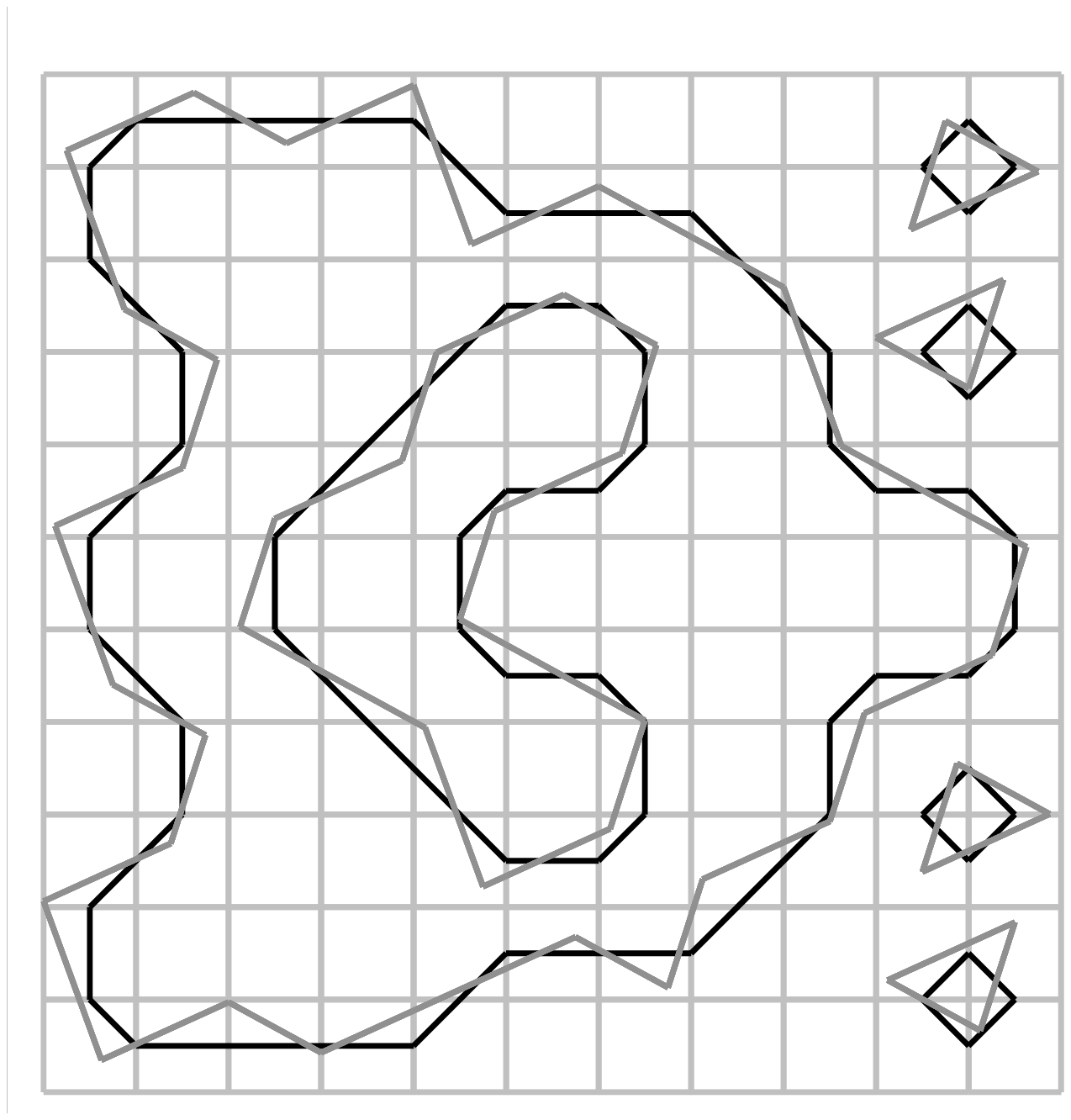}}
\newline
{\bf Figure 1.3:\/} The Quasi-Isomorphism Conjecture
in action for $3/8$.
\end{center}

For the parameter $3/8$,
Figure 1.3 superimposes the plaid polygons
in a certain region of the plane on top
of the said affine images of the
arithmetic graphs.

All this terminology will be defined in \S 9.
I believe that a generalization of my
proof of the Hexagrid Theorem in [{\bf S1\/}] will
establish the Quasi-Isomorphism Conjecture,
but the proof will require a separate paper.
The Quasi-Isomorphism Conjecture (once
proved) converts the plaid model
into a machine for generating results
about the arithmetic graphs associated
to outer billiards on kites.

The paper comes with a companion computer
program which illustrates many of the
results in this paper - in particular the
PET Equivalence Theorem.  One can download
this program from my website. The URL is \newline \newline
{\bf http://www.math.brown.edu/$\sim$res/Java/PLAID.tar\/}
\newline \newline
I discovered
all the results in the paper using a more
sophisticated version of the program above.
The more sophisticated version
is too complicated for public consumption.

I would like to think Peter Doyle, Pat Hooper, 
Sergei Tabachnikov, and
Ren Yi for a number of conversations about things
related to the plaid model.

\newpage

\section{The Grid Description of the Plaid Model}

\subsection{Basic Definitions}

\noindent
{\bf Even Rational Parameters:\/}
We will work with 
$p/q \in (0,1)$ with $pq$ even.
We call such numbers {\it even rational parameters\/}.
We will work closely with the auxiliary rationals
\begin{equation}
P=\frac{2p}{p+q}, \hskip 30 pt
Q=\frac{2q}{p+q}.
\end{equation}
These rationals have the property that $P+Q=2$ and $P/Q=p/q$.
\newline
\newline
{\bf Four Families of Lines\/}
We consider $4$ infinite families of lines.
\begin{itemize}
\item $\cal H$ consists of horizontal lines having integer $y$-coordinate.
\item $\cal V$ consists of vertical lines having integer $x$-coordinate.
\item $\cal P$ is the set of lines of slope $-P$ having
integer $y$-intersept.
\item $\cal Q$ is the set of lines of slope $-Q$ having
integer $y$-intersept.
\end{itemize}
Figure 2.1 shows these lines inside $[0,7]^2$
for $p/q=2/5$.  In this case, $P=4/7$ and $Q=10/7$.

\begin{center}
\resizebox{!}{3in}{\includegraphics{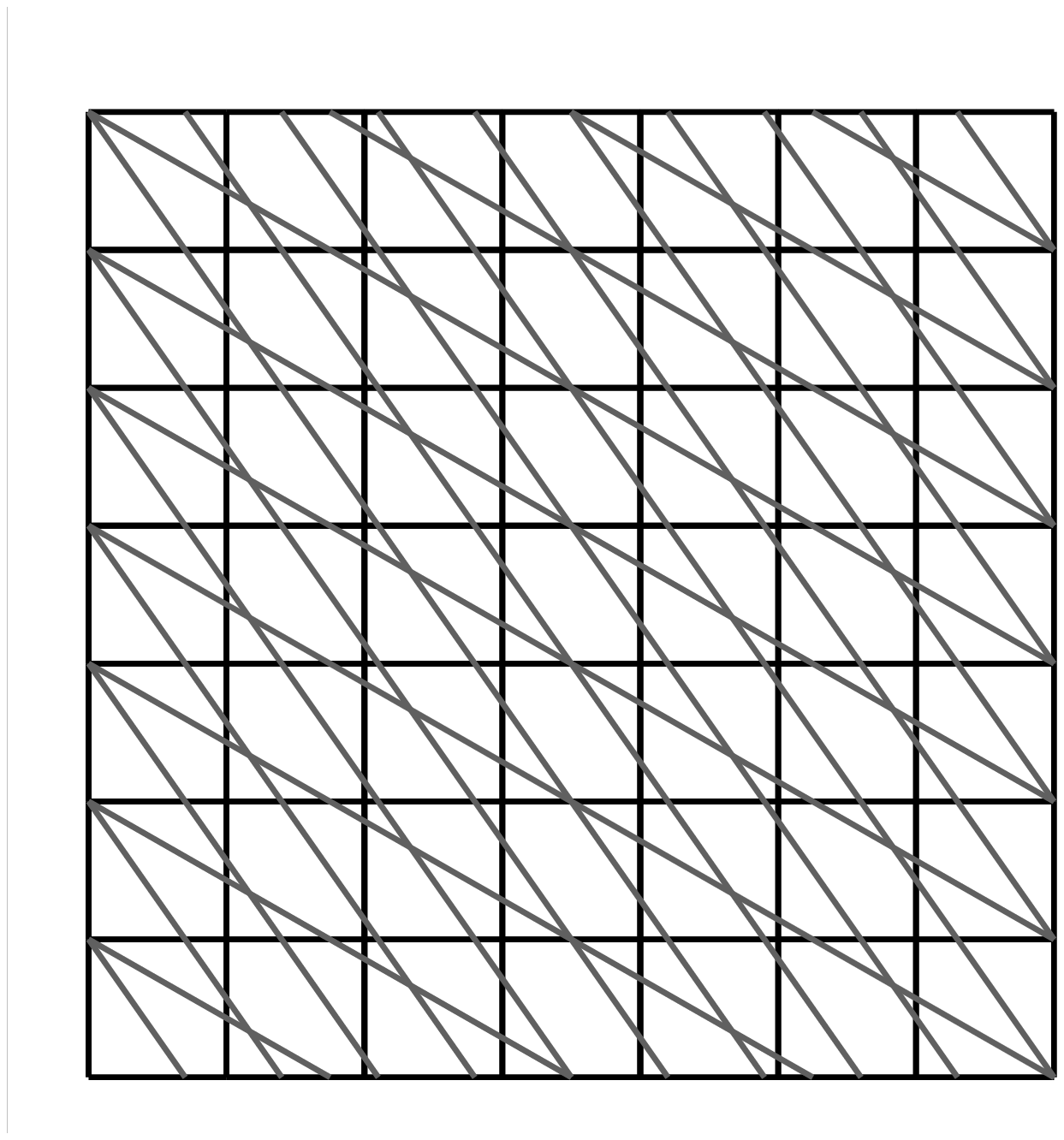}}
\newline
{\bf Figure 2.1:\/} The $4$ line families for $p/q=2/5$.
\end{center}

\noindent
{\bf Adapted Functions:\/}
Let $\Z_0$ and $\Z_1$ denote the
sets of even and odd integers respectively.
We define the following $4$ functions:
\begin{itemize}
\item $F_H(x,y)=2Py$ mod $2\Z$
\item $F_V(x,y)=2Px$ mod $2\Z$.
\item $F_P(x,y)=Py+P^2x+1$ mod 2$\Z$
\item $F_Q(x,y)=Py+PQx+1$ mod 2$\Z$.
\end{itemize}
We call $F_H$ and $F_V$ {\it capacity functions\/}
and $F_P$ and $F_Q$ {\it mass functions\/}.
We shall never be interested in the
inverse images of $\Z_1$, so
we will always normalize so that our
functions take values in $(-1,1)$.
Lemma \ref{anchor} gives a geometric
view of these functions, in terms of
circle rotations.
\newline
\newline
{\bf Intersection Points:\/}
We set
\begin{equation}
\omega=p+q, \hskip 30 pt
\Omega_j=\Z_j \omega^{-1}.
\end{equation}
We observe that, for each
index $A \in \{H,V\}$ and
$B \in \{P,Q\}$, we have
$F_A^{-1}(\Omega_0)=\cal A$ and
$F_B^{-1}(\Omega_1)=\cal B$.

We call $z$ an {\it intersection point\/}
if $z$ lies on both an $A$ line and a $B$ line.
In this case,
we call $z$ {\it light\/} if
\begin{equation}
\label{LIGHT}
|F_B(z)|<|F_A(z)|, \hskip 30 pt
F_A(z)F_B(z)>0.
\end{equation}
Otherwise we call $z$ {\it dark\/}.
We say that $B$ is the {\it type\/}
of $z$.  Thus $z$ can be light or dark,
and can have type P or type Q.  When
$z$ lies on both a $\cal P$ line and
a $\cal Q$ line, we say that $z$
has both types.
\newline
\newline
{\bf Capacity and Mass:\/}
For any index $C \in \{V,H,P,Q\}$, we
assign two invariants to a line
$L$ in $\cal C$, namely
$|(p+q)F_C(L)|$ and the sign of
$F_C(L)$.  We call the second invariant
the {\it sign\/} of $L$ in all cases.
 When $A\in\{V,H\}$, we call
the first invariant the {\it capacity\/} of $L$.
When $A \in \{P,Q\}$ we call the
first invariant the {\it mass\/} of $L$.
The masses are all odd integers in
$[1,p+q]$ and the capacities are
all even integers in $[0,p+q-1]$.
Thus, our rule can be described
as follows.
An intersection point is light
if and only if the intersecting
lines have the same sign, and the
mass of the one line is less than
the capacity of the other.
\newline

Figure 2.2 shows the light
intersection points inside the
suequare $[0,7]^2$ for the
parmeter $2/5$.  For each unit square
$S$ in this region, we have connected
the center of $S$ to the light points
on $\partial S$. 

\begin{center}
\resizebox{!}{4in}{\includegraphics{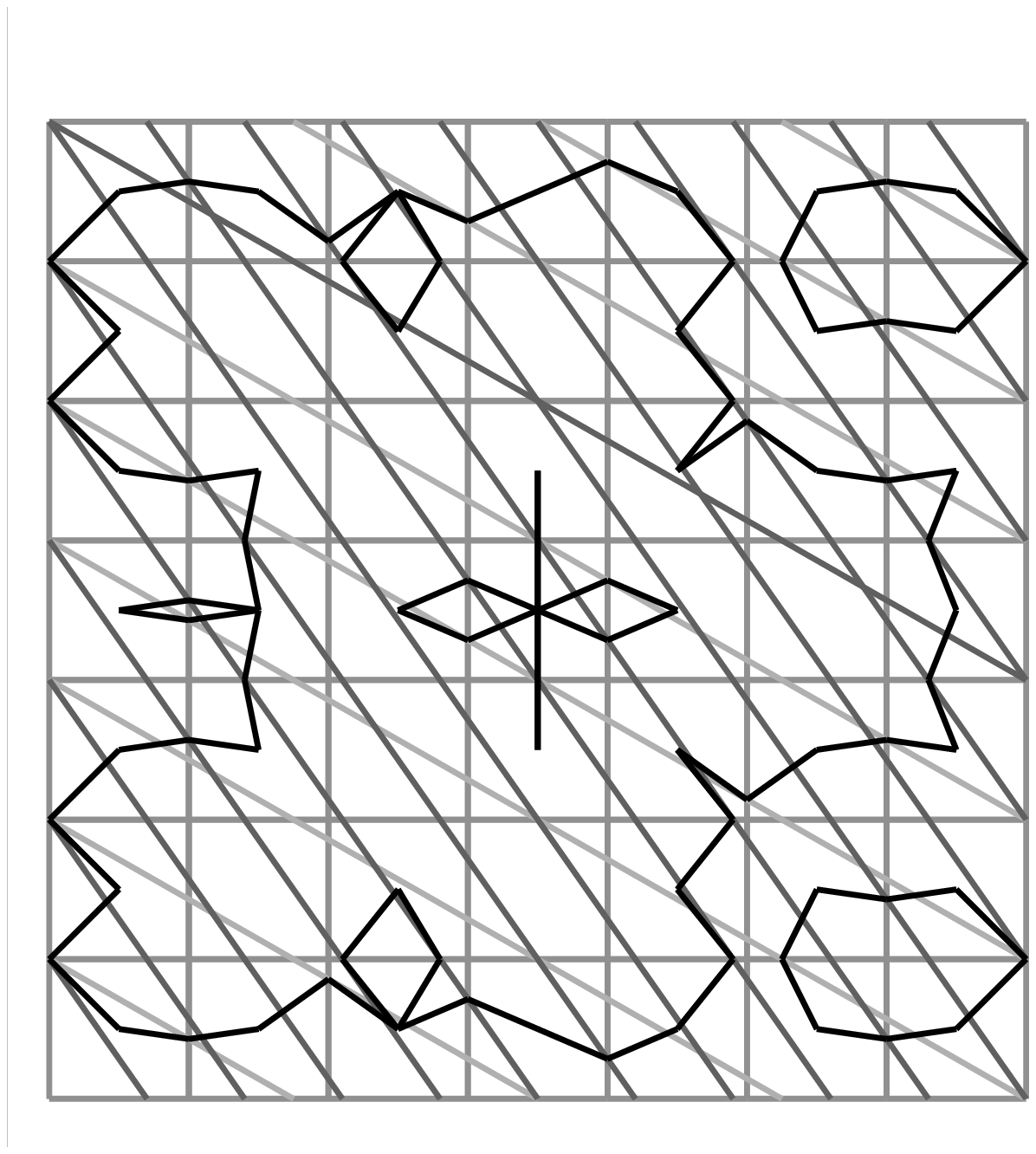}}
\newline
{\bf Figure 2.2:\/} Some light points for $p/q=2/5$.
\end{center}

 Notice that some
interesting curves seem to emerge.
Notice also that there seems to be a
small amount of junk  (in the form of
little loops) hanging off these curves.
Our rule below will prune away the
junk and keep the interesting part.
\newline
\newline
{\bf Double Counting Midpoints:\/}
We introduce the technical rule that
we count a light point twice if it
appears as the midpoint of a horizontal
unit segment.  The justification for this
convention is that such a point is 
always a triple intersection of a
$\cal P$ line, a $\cal Q$ line, and
a $\cal H$ line; and moreover this point
would be considered light either when
computed either with respect to the
$\cal P$ line or with respect to the
$\cal Q$ line.  See Lemma \ref{weak} below.
\newline
\newline
{\bf Good Segments:\/}
We call the edges of the unit squares
{\it unit segments\/}.  Such segments,
of course, either lie on $\cal H$
lines or on $\cal V$ lines.  We say
that a unit segment is {\it good\/} if it
contains exactly one light point.
  We say that a unit square
is {\it coherent\/} if it contains either
$0$ or $2$ good segments.  We say that the
plaid model is {\it coherent\/} at if all
squares are coherent for all parameters.
Here is the fundamental theorem concerning
the plaid model.

\begin{theorem}[Fundamental]
\label{fun1}
The plaid model is coherent for all parameters.
\end{theorem}

\noindent
{\bf Plaid Polygons:\/}
The Fundamental Theorem allows us to create a
union of embedded polygons.  In each
unit square we draw the line segment which connects
the center of the square with the centers of
its good edges.   The squares with no good edges
simply remain empty. We call these polygons
the {\it plaid polygons\/}. Figure 2.3 shows the
plaid polygons contained in $[0,7]^2$ for
the parameter $2/5$.

\begin{center}
\resizebox{!}{3.5in}{\includegraphics{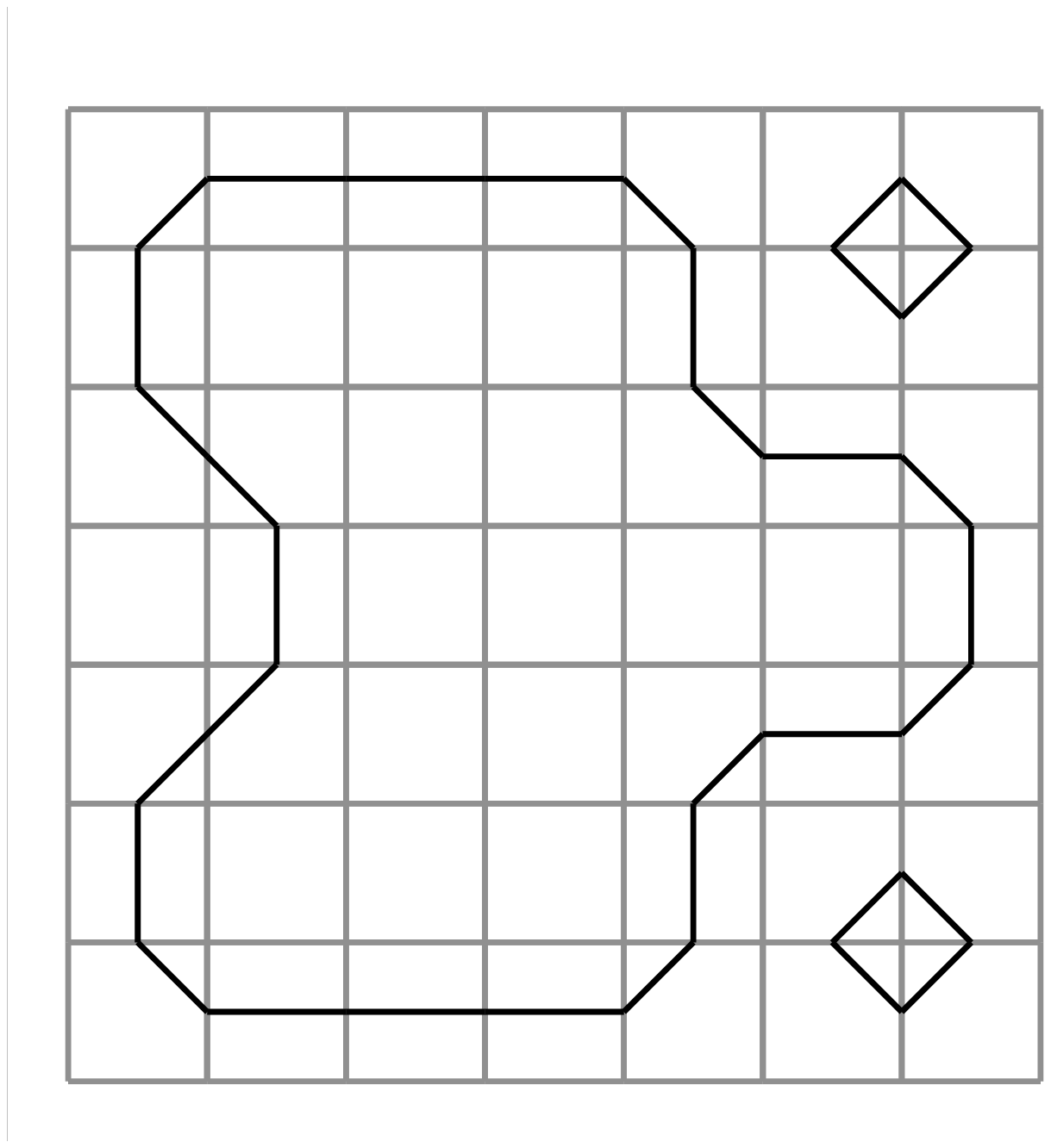}}
\newline
{\bf Figure 2.3:\/} The plaid polygons inside
$[0,7]^2$ for $p/q=2/5$.
\end{center}

\subsection{The Tune of the Model}
\label{tune}

Given a parameter $p/q$ let $\omega=p+q$.
We solve the equation
\begin{equation}
\label{tune1}
2 \alpha p \equiv \pm 1 \hskip 5 pt
{\rm mod\/} \hskip 5 pt \omega.
\end{equation}
and let $\alpha$ be which ever solution lies in
$(0,\omega/2)$.  We then define the {\it tune\/}
of the model to be
\begin{equation}
\tau(p/q)=\alpha/\omega.
\end{equation}
The numbers $p/q$ and $\tau(p/q)$ control
the geometry of the model, thanks to the
following lemma.

\begin{lemma}
\label{anchor}
For $k=0,...,(\omega-1)/2$, the
lines of capacity $2k$ have the form
$$x=k\tau \omega, \hskip 30 pt
x=\omega-k\tau \omega, \hskip 30 pt
y=k\tau \omega, \hskip 30 pt
y=\omega-k\tau \omega.$$
For $k=1,3,...,(\omega-1)$, the
the lines of mass $k$ have $y$-intercepts
$$
(0,k\tau \omega), \hskip 30 pt
(0,\omega-k\tau \omega).
$$
These equations are taken mod $\omega$.
\end{lemma}

\startproof
We will deal with the line $y=k\tau\omega$ and
the case $2p \alpha \equiv 1$ mod $\omega$.
The other cases are similar.  We compute
$$F_H(0,y)=[2Pk\tau \omega]_2=[4pk\tau]_2=
\frac{2}{\omega}[2pk\alpha]_{\omega}=
\frac{2}{\omega}[k]_{\omega}=\frac{2k}{\omega}.$$
We deal with the
$\cal P$ lines and $\cal Q$ lines at the same time.
We will deal with the case when the $y$-intercept
is $(0,ka)$.  The other case is similar.
We now are assuming that $k$ is odd.
We know that $[2Pky]_2=2k/\omega$.  We
Multiply this equation above by $\omega/2$, we see that
$[2pky]_{\omega}=[\omega Pky]_{\omega}=k.$
In short $2kpy \equiv k$ mod $\omega$.
Since $\omega$ and $k$ are odd, and
$2pky$ is even, we have
$2pky \equiv \omega+k$ mod $2\omega$.
But then
$$F_P(0,y)=[Py+1]_2=\frac{[\omega Pky+\omega]_{2\omega}}{\omega}=
\frac{[2pyk+\omega]_{2\omega}}{\omega}=\frac{k}{\omega}.$$
Hence the $\cal P$ and $\cal Q$ lines through
$(0,y)$ have mass $k$.
\endproof

The close connection between the plaid model and
circle rotations suggests that there ought to be
a lot of renormalization going on in the model.
We will not pursue this here, though we will exploit
the circle rotation property when we prove 
Theorem \ref{bigg2}.

\subsection{Symmetries}
\label{plaidsymm}

\noindent
{\bf The Symmtry Lattice\/}
We fix some even rational parameter $p/q$.
Let $\omega=p+q$ as above.
Let $L \subset \Z^2$ denote the
lattice generated by the two vectors
\begin{equation}
(\omega^2,0), \hskip 30 pt
(0,\omega).
\end{equation}
We call $L$ the {\it symmetry lattice\/}.
\newline
\newline
{\bf Blocks:\/}
We define the square $[0,\omega]^2$
to be the {\it first block\/}.
The pictures above always show the first block.
In general, we define a {\it block\/} to be
a set of the form $B_0 + \ell$, where $B_0$ is the
first block and $\ell \in L$.
With this definition, the lattice $L$
permutes the blocks.  
We define the {\it fundamental blocks\/} to
be $B_0,...,B_{\omega-1}$, where
$B_0$ is the first block and
\begin{equation}
B_k=B_0+(k\omega,0).
\end{equation}
The union of the fundamental blocks is a
fundamental domain for the action of $L$.
We call this union the {\it fundamental domain\/}.
\newline
\newline
{\bf Translation Symmetry:\/}
All our functions $F$ are $L$-invariant.
That is, $F(v+\ell)=F(v)$ for all $\ell \in L$.
This follows immediately from the definitions.
Note that a unit segment in the
boundary of a block does not have any
light points associated to it
because $F_H=0$ on the horizontal
edges of the boundary and $F_V=0$ on
the vertical edges of the boundary. 
(A vertex of such a segment might be a
light point associated to another
unit segment incident to it but this
doesn't bother us.)
Thus, assuming the truth of the
Fundamental Theorem -- so that the definition
of a plaid polygon makes sense -- every plaid polygon is
contained in a block, and is translation equivalent
to one that is contained in a fundamental block.
\newline
\newline
{\bf Rotational Symmetry:\/}
The plaid model is invaruant under
reflection in the origin. 
Here is how we see this.
For each of the functions $F$, we have
$F(-x,-y)=-F(x,y)$ mod $2\Z$.  This
means that the rotation in the origin
preserves the set of light 
intersection points and the set of
dark intersection points, and preserves
the types.
\newline
\newline
{\bf Reflection Symmetry:\/}
Let $\widehat L$ denote the group
of isometries generated by
reflections in the horizontal and
vertical midlines (i.e. bisectors)
of the blocks.   We will show
through a series of lemmas that
the elements of $\widehat L$ 
preserve the set of plaid polygons.
However, for the reflections in 
horizontal lines, the vertical
interection points of type P are
swapped with the vertical intersection
points of type Q.

The group $\widehat L$ is generated
by $L$, by the rotation mentioned above,
and by reflection in the
$x$-axis. So, to finish the proof, we
just have to analyze what happens for
reflection in the $x$-axis.  This
is the interesting case.

\begin{lemma}
\label{REFLECT}
The set of light points for the plaid
model is invariant under reflection
in the $x$-axis.  This reflection
preserves the type of the horizontal
light points and swaps the types of
the vertical light points.
\end{lemma}

\startproof
We will show that the reflection
preserves a light point
$z=(x,y)$ of type P.
The proofs in the other cases
are very similar.

Suppose first that $z$ lies on a vertical edge.
Let $z'=(x,-y)$.  We have
\begin{equation}
F_V(z')=F_V(z), \hskip 15 pt
F_P(z')=F_V(z)-F_Q(z) \hskip 15 pt
F_Q(z')=F_V(z)-F_P(z).
\end{equation}
The last two equations hold mod $2\Z$.
Using the rotation symmetry, we reduce to the
case when $F_P(z)>0$.  But then $F_V(z)>0$ and
$F_P(z')+F_Q(z)=F_V(z)$ when the values are
forced to lie in $(-1,1)$.
Note that $\Z_0-\Z_1=\Z_1$.  So,
$F_Q(z') \in \Z_1$.  But then our
definitions tell us that $z'$ is a
light point of type Q.

Now suppose that $z$ lies on a horizontal segment.
We have
\begin{equation}
\label{reflect2}
F_H(z')=-F_H(z), \hskip 10 pt
F_P(z')-F_H(z')=F_P(z), \hskip 10 pt
F_Q(z')-F_H(z')=F_Q(z).
\end{equation}
Using the rotational symmetry, we can assume that
\begin{equation}
\label{reflect3}
F_H(z)<F_P(z)<0
\end{equation}
   But then
$F_H(z')>0$.  Equation \ref{reflect2} forces
$F_P(z')<F_H(z')$.  Equation \ref{reflect3} 
combines with the equation
$$F_P(z')=F_P(z)-F_H(z) \hskip 30 pt {\rm mod\/}\ 2\Z$$
to force $F_P(z')>0$. This
$0<F_P(z')<F_H(z')$.  Finally, Equation
\ref{reflect2} tells us that $F_P(z') \in \Omega_1$.
So, $z'$ is a light point of type P.
\endproof

\subsection{The Number of Intersection Points}

The purpose of this section is to prove
the following result.

\begin{lemma}
\label{2intersect}
Each unit segment contains $2$ intersection points.
\end{lemma}

\startproof
Let $e$ be a vertical edge.  Let $L$ be the
vertical line through $e$.
There is some $\alpha \in [0,1]$ such that
the intersection points of type $P$ along
$L$ have the form $n+\alpha$, where
$n \in \Z$.  The same goes for the 
points of type Q.  Hence, there are
exactly two of them in $e$.  (This works
even when $\alpha=0$.)

Now let $e$ be a horizontal edge. 
If we forget about whether the intersection
points are light or dark, the whole picture
is symmetric under translation by $(0,1)$
and also $(p+q,0)$.  
So, we can assume that $e$ lies on the
south boundary the first block.
The intersection points of type $P$ have the
form $(n/P,0)$ where $n\in \Z$ and the intersection
points of type $Q$ have the form $(n/Q,0)$.
\newline
\newline
{\bf Case 1:\/}
If $e$ is the central edge, then $e$
contains the intersection points
$$(p/P,0)=(q/Q,0)=((p+q)/2,0).$$
This common point is counted twice, by convention.
\newline
\newline
{\bf Case 2:\/}
If $e$ is the westernmost edge, then
$e$ contains the two intersection points
$(0,0)$ and $(1/Q,0)$.  If $e$ is the
easternmost edge, then $e$ contains
$(p+q,0)$ and $(p+q,0)-(1/Q,0)$.
\newline
\newline
{\bf Case 3:\/}
If $e$ is not one of the edges above,
then neither the boundary nor the
midpoint of $e$ contains an
intersection point.  Since $1/Q \in (1,2)$
we know that $e$ contains at least $1$
point of type Q and at most $2$ of them.
We will show that if $e$ does not contain
a second point of type Q then $e$ contains
a point of type P.
Let $(k_1/Q,0)$ be the point of
type Q that $e$ does contain.
We must have
$k_1 \in (Qm+Q-1,Qm+1)$,
for otherwise we could add or
subtract $1$ from $k_1$ and
produce another intersection
point of type Q in $e$.
We claim that there is some 
point of type P inside $e$.
We seek a point 
$k_2 \in (Pm,Pm+P)$.
This time we have
$$k_1+k_2 \in (2m+Q-1,2m+P+1)=(2m+1-P,2m+1+P),$$
The value $k_2=(2m+1)-k_1$ does the job.
\endproof

\subsection{Capacity and Mass}

Now we come to a more subtle result
which suggests the hierarchical
nature of the plaid model.  The lines of small
capacity have very few light points, so they
predict something about the large scale
geometry of the loops in the model.  As we add
more lines of higher capacity, the picture
of the loops fills in at finer scales.
We will take up this discussion in detail in
\S \ref{loops}.

\begin{theorem}
\label{hier}
Let $B$ be any block.
For each even $k \in [0,p+q]$ there are
$2$ lines in $\cal H$ and $2$ lines
in $\cal V$ which have capacity $k$ and
intersect $B$.  Each such line
carries $k$ light points in $B$.
\end{theorem}

\begin{lemma}
Statement 1 of Theorem \ref{hier} is true.
\end{lemma}

\startproof
Recall that $\omega=p+q$.
Given the periodicity of the functions $F_H$ and
$F_V$, it suffices to prove this for the
first block.  We will prove the result for
$\cal H$. The result for $\cal V$ has virtually
the same proof.
We are simply trying to show that there are
exactly $2$ integer values of $y$ in $[0,\omega]$
such that $4py=\pm k$ mod $2\omega$.
Writing $k=2h$, we see that this equation
is equivalent to $2py = \pm h$ mod $\omega$.
This has $2$ solutions
mod $\omega$ because $2p$ is relatively prime
to $\omega$.
\endproof

\begin{lemma}
\label{strong}
Theorem \ref{hier} holds in the vertical case.
\end{lemma}

\startproof
Let $L$ be a vertical line of capacity $k$.
The case $k=0$ is trivial, so we assume $k>0$.
In this case, no point of type P coincides
with a point of type Q.  We will show that
there are $k/2$ light points of type P in
$L \cap B$.  By reflection symmetry, the same
goes for the points of type Q.

Let $S_P \subset \Z$ denote those points $m$ such
that the line $\Lambda_m$ of $\cal P$ through $(0,m)$ 
intersects $L \cap B$.
The set $S_P$ is the intersection of
$\Z$ with a segment of length $\omega$.  When $k>0$ this
segment has endpoints which are not in $\Z$. 
Hence $S_P$ consists of exactly $\omega$ consecutive
integers.

The restriction
of $\omega F_P$ to $\Lambda_m$ is $2pm+\omega$.
Choose any odd $\ell \in (0,k)$.  To prove our
result we just have to show that there is exactly
one $m \in S_P$ such that 
$$2pm+\omega \equiv \ell \hskip 20 pt
{\rm mod\/}\ 2\omega$$
Since $\omega$ is odd, we can write $\ell=2d+\omega$, and the
congruence above is equivalent to
$pm \equiv d$ mod $\omega$.  There is exactly one solution
to this in any run of $\omega$ consecutive integers.
\endproof

\begin{lemma}
\label{weak}
Theorem \ref{hier} holds in the horizontal case.
\end{lemma}

\startproof
Let $L$ be a horizontal line of capacity $k$
which interects the block $B$.
By symmetry,
we can assume that $B$ is one
of the fundamental blocks.
Let $L_+$ and $L_-$ be the two horizontal lines
of capacity $k$ and $-k$ respectively.
Reflection in the horizontal midline of $B$
swaps these two lines.

Let $\ell \in (0,k)$ be some odd integer.
Let $S(\ell)$ denote the set $x \in \R$ with the
following property.  There is some $y \in \R$
such that either
\begin{itemize}
\item $(x,y)$ is a light point of type P on $L_-$.
\item $(x,y)$ is a light point of type Q on $L_+$.
\end{itemize}
We will show that $S(\ell)$ has cardinality $2$
for all $\ell=1,...,k-1$.
This fact, together with the bilateral symmetry,
establishes the lemma.

Let $S_Q$ (respectively $\widehat S_P$) denote the set
$m \in \Z$ such that the line of slope $-Q$
(respectively $+P$) through $(0,m)$ intersects $L_+$.
By symmetry, our set $S(\ell)$ has the following
description.  We consider all points 
$m \in S_Q \cup \widehat S_P$ having
$\omega F_P(0,m) \equiv \ell$ mod $2\omega$.
We then intersect the appropriately sloped
line through $(0,m)$ with $L_+$ and
record the $x$ coordinate.  We call
$(x,y_+)$ a {\it guide point\/}.
Here $y_+$ is the $y$-coordinate of $L_+$.

We claim $\widehat S_P \cup S_Q$ contains
numbers in every congruence class mod $2\omega$.
Consider first the case when $B$ is the first block.
In this case, $y_+$ is both the right endpoint of
$\widehat S_P$ and the left endpoint of $S_Q$.   So,
in this case, the union is a run of $2\omega+1$ consecutive
integers. So, the congruence property holds.
When we replace the first block $B$ by the
$k$th block $B'$, the set $\widehat S_P$ moves down
by $2pk$ units and the set $S_Q$ moves up by $2qk$
units.  Hence, the congruence property still holds.

As in the proof of Lemma \ref{weak}, we are trying
to solve the equation
$$pm \equiv d\hskip 5 pt {\rm mod\/}\hskip 5 pt \omega,
\hskip 30 pt \ell=2d+\omega, \hskip 30 pt
m \in \widehat S_P \cup S_Q.$$
The congruence property implies that
there are $2$ solutions $m$ and $m'$
modulo $2\omega$.

If $m \equiv m'$ mod $2\omega$, then the
corresponding giode points coincide
because $P+Q=2$ and relevant lines of slope
$-Q$ and $P$ intersect the $y$-axis at
points which are $2k\omega$ apart.
If $m \equiv m'+\omega$ mod $2\omega$ and
$y_+=(m+m')/2$ the two guide points lie at the
center point of $L_+ \cap B$ and our convention
says to count the point twice.
In all other cases, the two guide
points are distinct.
\endproof

\subsection{Remote Adjacency and Particles}
\label{particle}

The results in the previous section are
a kind of conservation principle.  As
we move from block to block, the number of
light points on a line of capacity $k$
(namely $k$ of them) does not change.
In this section, we further the physical
analogy and explain how to think about
our intersection points as moving particles.
This is a key step in our proof of
the Fundamental Theorem.   

Let $B_0,...,B_{\omega-1}$ be the
fundamental blocks.
Let $a \in (0,\omega)$ be such that
\begin{equation}
2ap \equiv -1 \hskip 20 pt
{\rm mod\/}\ \omega.
\end{equation}
We say that the blocks
$B_j$ and $B_{j+a}$ are {\it remotely adjacent\/}.
These indices are taken mod $\omega$.
We write $B_j \to B_{j+a}$.  Note that
$a$ is relatively prime to $\omega$, so that
cycle $B_0 \to B_a \to B_{2a}...$ lists
out every fundamental block.

Given a fundamental block $B$, there is a 
horizontal translation
$T$ such that $T(B)=B_0$.   Given any point
$z \in B$, we define $[z]=T(z) \in B_0$.
The point $[z]$ records the
position of $z$ within $B$.

\begin{lemma}
Suppose $B \to B'$.
Let $H$ be a $\cal H$ line which intersects $B$ and $B'$.
Let $z$ be an intersection point on $H$.
Suppose $z$ has type P (respectively type Q) and is not
on the left (respectively right) edge of $B$.
Then there is a point $z' \in B' \cap H$ of the same type and
brightness such that
$[z']-[z]=(P^{-1},0)$ (respectively $[z']-[z]=(-Q^{-1},0)$.)
\end{lemma}

\startproof
Consider first the case when $z$ has type P.
The first thing we have to do is check that
$z' \in B'$. 
The lines in $\cal P$ intersect $H$ in
points of the form $nP^{-1}$ for $n \in \Z$.
So, $z$ must be at least $P^{-1}$ from
the right edge of $B$.  This gives $z' \in B'$.

Now we will show that $z'$ is light
if and only if $z$ is light.
Let $L$ and $L'$ respectively be the
$\cal P$ lines through $z$ and $z'$
respectively. Let $y$ and $y'$ be
such that $(0,y) \in L$ and $(0,y') \in L'$.
The blocks $B$ and $B'$ differ by
a horizontal translation of $a(\omega)$.
Given that the lines in $\cal P$ have
slope $-P=-2p/(\omega)$, we get
$y'-y=2pa+1$.  But this difference
is $0$ mod $\omega$, so that
function $F_P$ gives the same value
to points on $L$ and points on $L'$.
Now our basic criterion says that
$z'$ is light if and only if $z$
is light.

When $z$ has type Q, the proof is the same
except that $y'-y=2qa-1$
and $2aq \equiv -1$ mod $\omega$. 
\endproof

In the vertical case,
the situation is different
because each vertical line intersects
at most one fundamental block.
So, we consider a family of vertical
lines $\{V_k\}$ such that
$V_k$ intersects $B_k$ in the same
relative position for all
$k=0,...,\omega-1$.  We could say that
$[B_k]$ is the same line, independent
of $k$.  We write $V_k \to V_{k+a}$.
For the next result, it is useful to
think of our blocks as cylinders,
with the tops and bottoms identified.

\begin{lemma}
Let $B \to B'$ be fundamental blocks.
Let $V \to V'$ be vertical lines which
respectively intersect these blocks.
Let $z$ be an intersection point of type P
(respectively type Q) on $V$.
Then there is an intersection point $z' \in V'$ of
the same type and brightness so that
$[z']=[z]+(0,1)$ (respectively
$[z']=[z]-(0,1)$.)
\end{lemma}

\startproof
This has the same proof as in the horizontal case.
\endproof

\noindent
{\bf Horizontal Motion:\/}
When $(z,z')$ satisfy the relation
indicated in one of the two lemmas above,
we write $z \to z'$.
Consider the horizontal case first.
Let $z_0$ be some intersection point
which starts out in the left edge of a fundamental block.
We write $z_0 \to ... \to z_{\omega}$.
The points $z_0,...,z_{2p-1}$ all have type P,
and their translates
$[z_0],...,[z_{2p}]$ move east across the
first block by in steps of $P^{-1}$.
then $z_{2p}$ lies in the right edge of
its block and has type Q.  The
points $z_{2p},...,z_{\omega-1}$ all have
type Q and their translates
$[z_{2p}],...,[z_{\omega-1}]$ move west across the
first block by in steps of $P^{-1}$.
After $\omega$ steps the cycle is done.
\newline
\newline
{\bf Vertical Motion:\/}
We again have the cycle
$z_0,...,z_{\omega-1}$.  This time the
type does not change. The translates
$[z_0],...,[z_{\omega-1}]$ move
north in steps of $1$ unit in the type P
case, and south in steps of $1$ unit
in the type Q case.
\newline
\newline
{\bf Particles:\/}
We call the points $z_0,...,z_{\omega-1}$
comprising these cycles {\it particles\/}.
We call each of the individual points
$z_i$ {\it instances\/} of the particle.
So, a particle consists of all its instances.
Thus, a particle consists of $\omega$ intersection
points.  These points are all either light or dark,
and so we can speak of {\it light particles\/}
and {\it dark particles\/}.
The horizontal particles have points of both
types, depending on their direction of motion,
and the vertical particles have points all
of the same type.
If $\omega$ is large and we rescale the picture
so that the individual blocks have unit size,
the ``movement'' of the particles as we cycle
through the blocks looks very much like 
locally linear motion.

\newpage

\section{The Tile Description of the Plaid Model}

\subsection{Classifying Pairs}

\noindent
{\bf Tiles and Connectors:\/}
Let $Q$ be a unit square in the
usual square grid.  We label the
edges of $Q$, in the obvious way,
with the letters N(orth), S(outh), E(ast), and W(est).
 There are $7$ ways to
draw an edge from the midpoint of one edge of $Q$ to
the midpoint of another edge
One of the ways is that we simply draw no edge
at all, and the other $6$ ways correspond to 
unordered pairs in the set
$\{N,S,E,W\}$.   By way of example, we
say that the {\it NW-tile\/} is the one which has
an edge joining the N edge to the W edge.
We say that an {\it empty tile\/} is a tile
with no edge drawn.  In the $6$ cases when we
actually draw something, we call this segment a
{\it connector\/}. 
\newline
\newline
{\bf Coherent Tilings:\/}
Suppose that we have two adjacent tiles sharing
a common edge.  We say that these tiles
{\it match\/} if the common edge is involved
in the connectors of both tiles, or in neither.
Suppose we have a tiling of the plane which
uses the various tiles.  We call the tiling
{\it coherent\/} if every pair of adjacent
tiles match across their common edge.
\newline
\newline
{\bf Classifying Spaces:\/}
We say that a {\it classifying space\/}
is an convex polytope which has been partitioned into a
finite number of smaller convex polytopes, each of which
has been given one of $7$ labels corresponding
to the different types of tiles.
By {\it partition\/} we mean
that the polytopes have pairwise
disjoint interiors.  We call the classifying
space {\it integral\/} if every convex
polytope in sight is an integer convex polytope - i.e.,
the convex hull of a finite union of integer vectors.
\newline
\newline
{\bf Classifying Pairs:\/}
Recall that $\Z_1$ is the set of odd integers.
Let $\cal C$ denote the set of centers of
squares in the usual square grid.  The points
of $\cal C$ have the form $(a/2,b/2)$ where
$a,b \in \Z_1$.  Let $X$ be a classifying
space and let
$\Xi: \R^2 \to X$ be a map.  We call
$\Xi$ a {\it classifying map\/} if
\begin{itemize} 
\item $f$ is entirely defined on $\cal C$.
\item $f$ maps $\cal C$ into the union of
interiors of the pieces of the partition.
\end{itemize}
We say that $(\Xi,X)$ is a {\it classifying pair\/}.

\subsection{Description of the Space}
\label{classmap}

We fix $P \in [0,1]$ and let
$\Lambda_P$ denote the lattice generated by the
vectors
\begin{equation}
(0,2,P,P), \hskip 30 pt (0,0,2,0), \hskip 30 pt (0,0,0,2).
\end{equation}
The action on the first coordinate is trivial.

We take the quotient
\begin{equation}
X_P=(\{P\} \times \R^3)/\Lambda_P.
\end{equation}
The cube
$\{P\} \times [-1,1]^3$
serves as a fundamental domain for the action of
$\Lambda_P$ on $\{P\} \times \R^3$. However, the
boundary identifications depend on $P$.

We define
\begin{equation}
X=\bigcup_{P \in [0,1]} X_P.
\end{equation}
$X$ is a flat affine
manifold.
Let $\widehat \Lambda$ denote 
the abelian group of affine transformations
$\widehat \Lambda=\langle T_1,T_2,T_3 \rangle$, where
\begin{enumerate}
\item
$T_1(x_0,x_1,x_2,x_3)=(x_0,x_1+2,x_2+x_0,x_3+x_0)$.
\item
$T_2(x_0,x_1,x_2,x_3)=(x_0,x_1,x_2+2,x_3)$.
\item
$T_3(x_0,x_1,x_2,x_3)=(x_0,x_1,x_2,x_3+2)$.
\end{enumerate}
Then
\begin{equation}
X=([0,1] \times \R^3)/\widehat \Lambda.
\end{equation}
$[0,1] \times [-1,1]^3$ serves as a
fundamental domain for $X$.
\newline

We have a map
\begin{equation}
\Xi: [0,1] \times \R^2 \to X,
\end{equation}
defined by
\begin{equation}
\label{map}
\Xi(P,x,y) = (P,2Px+2y,2Px,2Px+2Py)
\hskip 10 pt {\rm mod\/}\hskip 5 pt \Lambda_P.
\end{equation}
When the choice of $P$ is understood, we will
think of $\Xi_P$ as a map
from $\R^2$ into $X_P$.  Sometimes
we will drop off the first coordinate, namely $P$,
and think of $X_P$ as a subset of $\R^3$.
For each even rational parameter $p/q$ we
set $P=2p/(p+q)$, as above, and
restrict $\Xi_P$ to the set
$\cal C$ of square tiles.  

\subsection{Checkerboards}
\label{checker}

Now that we have described our classifing map
and the classifying space, we need to
describe the partition of the classifying space.
As a prelude, we describe a certain family
of partitions of the square $[-1,1]^2$.
The decomposition is based on a triple
$U=(u_1,u_2,u_3)$ of
numbers
\begin{equation}
-1 \leq u_1 \leq u_2 \leq u_3 \leq 1
\end{equation}
and also on a $4 \times 4$ ``permutation matrix'' $M$ whose
nonzero entries are replaced by the symbols
N,S,E,W.  This ''matrix'' is just a combinatorial
object for us.  It is useful to think of the triple
of numbers as being pairwise distinct, and we will
draw pictures this way.  The cases when some of them
coincide, or equal $\pm 1$, are degenerate limits.

Forgetting about the symbols, our matrix will
always be symmetric with respect to the
diagonal line of slope $1$.  Also, there is a
relation between $M$ and $U$
that is best discussed after the construction.
Here is the $3$-step construction.

\begin{enumerate}
\item We partition $[-1,1]^2$ into a $4 \times 4$
checkerboard by inserting the lines
$x=u_1,u_2,u_3$ and also the lines
$y=u_1,u_2,u_3$.
\item There are $4$ squares of the grid which
correspond go the nonzero entries of $M$. We
label these according to the symbols in $M$.
We call these the special rectangles.
\item In the row and column of a given special
rectangle, we add the symbol corresponding to
that rectangle.  Thus, each of the other $12$
rectangles is labeled by an unordered pair of
symbols from $\{N,S,E,W\}$.  
\end{enumerate}

\noindent
{\bf Compatibility:\/}
$M$ and $U$ must be compatible
in the sense that the
$4$ special rectangles are squares.
\newline  

As an illustration, consider the input data
$$(u_1,u_2,u_3)=(-2/3,0,1/3), \hskip 30 pt
M=\left[\matrix{
0&0&0&W \cr
N&0&0&0 \cr
0&E&0&0 \cr
0&0&S&0}\right].
$$
For this choice of $M$, the
compatibility condition is that
$-u_1-u_2+u_3=1$.
Here is the partition we get.

\begin{center}
\resizebox{!}{3.5in}{\includegraphics{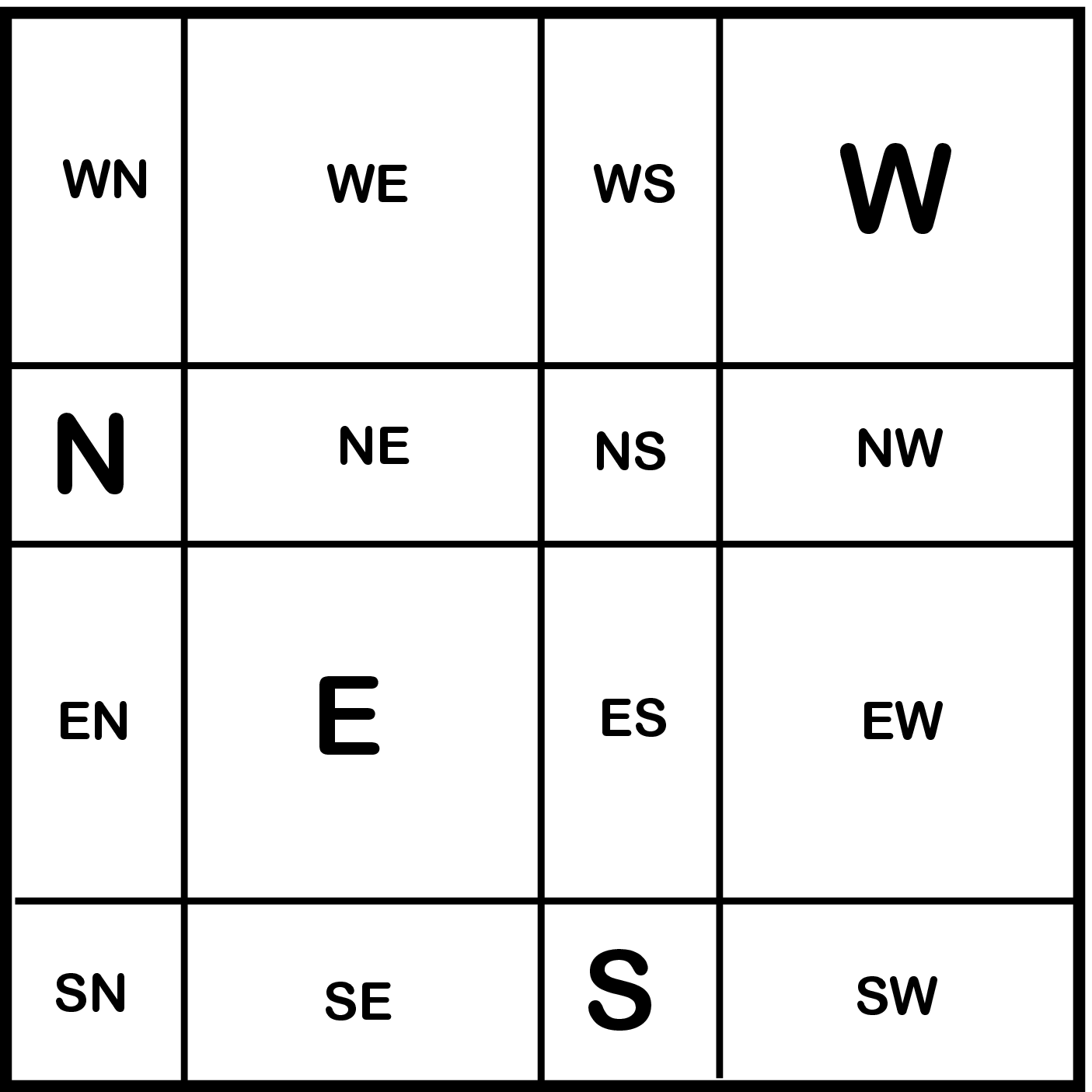}}
\newline
{\bf Figure 3.1:\/} The checkerboard partition
\end{center}

For each choice $M$, we can consider the set $U_M$ of all
triples which are compatible with $M$.  In the example
above, and in all the cases considered below,
$U_M$ is a triangle with integer vertices.
We think of the product 
\begin{equation}
U_M \times [-1,1]^2
\end{equation}
as a fiber
bundle over $U_M$. In each fiber we make the
checkerboard construction.  The union of all subsets
having the same label is an integer convex polytope,
and $U_M \times [-1,1]^2$ is also an integer convex
polytope.  We denote this space by $\langle M\rangle$

\subsection{Description of the Partition}
\label{partition0}

We will describe our partition of $X$ in
terms of subsets of the fundamental
domain $[0,1] \times [-1,1]^3$.  The
coordinates we use on $X$ are
$(P,T,U_1,U_2)$.

Let $X_{PT}$ denote the $2$ dimensional slice
obtained by holding $P$ and $T$ fixed. Our
partition is such that $X_{PT}$ has a checkerboard
partition, as discussed above.  Now we will
specify the parameters for the checkerboard
partition as functions of $P$ and $T$.

There are $3$ separate ranges,
so we will describe $3$ sets of data.

\begin{equation}
\label{zone1}
T \in [\!-\!1,\!-\!1\!+\!P]: \hskip 10 pt
U=(T,1\!-\!P,2\!-\!P\!+\!T), \hskip 10 pt
M=\left[\matrix{
0&0&0&W \cr
N&0&0&0 \cr
0&E&0&0 \cr
0&0&S&0}\right].
\end{equation}

\begin{equation}
\label{zone2}
T \in [\!-\!1\!+\!P,1\!-\!P]: \hskip 10 pt
U=(\!-\!1\!+\!P,T,1\!-\!P), \hskip 10 pt
M=\left[\matrix{
N&0&0&0 \cr
0&0&E&0 \cr
0&W&0&0 \cr
0&0&0&S}\right].
\end{equation}

\begin{equation}
\label{zone3}
T \in [1\!-\!P,1]: \hskip 10 pt
U=(\!-\!2\!+\!P\!+\!T,\!-\!1\!+\!P,T), \hskip 10 pt
M=\left[\matrix{
0&N&0&0 \cr
0&0&W&0 \cr
0&0&0&S \cr
E&0&0&0}\right].
\end{equation}

We can picture the data by fixing some
value of $P$ and plotting the graphs
of the functions $(u_1,u_2,u_3)$ as
a function of $T$.  Figure 3.2 shows the picture
for three values of $P$.
Here $Q$ is such that $P+Q=2$.

\begin{center}
\resizebox{!}{1.8in}{\includegraphics{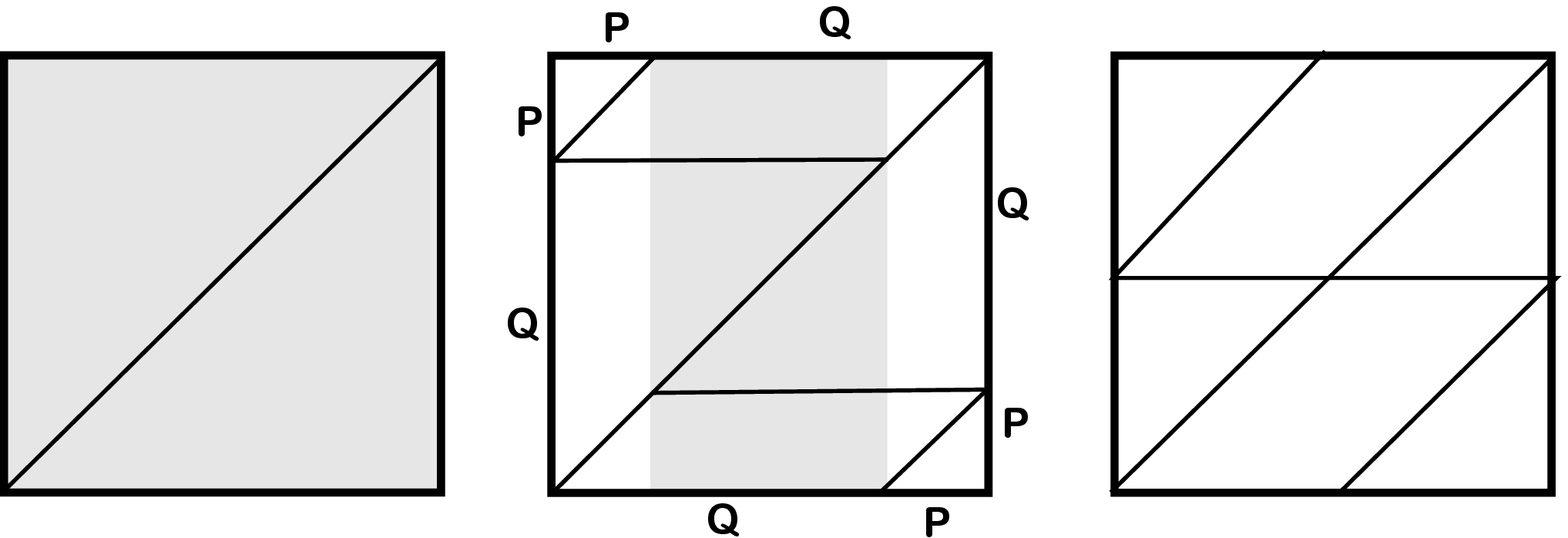}}
\newline
{\bf Figure 3.2:\/} The checkerboard data for
$P=0$ and $P=1/2$ and $P=1$.
\end{center}

The space $X$ is partitioned into $7$ regions.
$6$ of the pieces are obtained as follows:
We choose some label such as NW, and consider
the union, over all fibers, of the rectangles
having the label NW.  We call this piece
$X_{NW}$ and we give similar
names to the other pieces. What remains is $X_0$.

If we restrict $P$ to $(0,1)$, the
partition, combinatorially
speaking, has a product structure.  There are
no vertices at all.  When $P=0$ the only
vertices occur when $T \in \{-1,1\}$.  When
$P=1$ the only vertices occur when $T \in \{-1,0,1\}$.
Here is a list of all the vertices.

vertices
\begin{itemize}
\item $(0,a,b,c), \hskip 30 pt a,b,c \in \{-1,1\}$.
\item $(1,a,b,c), \hskip 30 pt a,b,c \in \{-1,0,1\}$.
\end{itemize}

By studying Figure 3.2 we can reconstruct each
individual polytope.  There are $16$ total, and
each one intersects the generic fiber in a
rectangle.  Indeed, each polytope is itself
the closure of a ``rectangle bundle'' over an open
triangle. 

Here is another way to think about it.
Recall that $X$ is segmented into $3$ ``zones'',
depending on the relevant values of $T$ with
respect to $P$, as in Equations
\ref{zone1} - \ref{zone3}.  Each zone has a matrix
$M$ attached it.  The linear isomorphism
\begin{equation}
(P,T,U_1,U_2) \to (U,U_1,U_2)
\end{equation}
maps each zone isomorphically to $\langle M \rangle$,
for the relevant choice of $M$.
Here $U$ is the triple in the relevant
choice of Equation \ref{zone1}, Equation
\ref{zone2} or Equation \ref{zone3}.

\subsection{Well Definedness}
\label{defined}

In this section, we explain why our classifying
pair $(\Xi,X)$ produces a well-defined tiling
for each even rational parameter $p/q$.  The issue
is that perhaps $\Xi$ maps some tile center
into the boundary of the partition.  We
will rule this out, and in fact will prove 
somewhat more about the map.

Let $\omega=p+q$ and $P=2p/\omega$.
We will usually drop the reference to the parameter
in our notation.  Also, we will forget the
$P$ coordinate and work in $\R^3$.  Recalling
that $\Z_0$ and $\Z_1$ respectively denote the
set of even and odd integers, we define
\begin{equation}
\label{points}
{\cal X\/}=\bigg(\frac{\Z_1}{\omega},\frac{\Z_0}{\omega},\frac{\Z_0}{\omega}
\bigg).
\end{equation}
Note that the lattice $\Lambda$ consists entirely
of vectors in
$(\Z_0/\omega)^3$.
Hence $\Lambda$ preserves $\cal X$.
Recall that $L$ is the lattice of symmetries
of the plaid model, and $\Xi$ is well
defined on the quotient ${\cal C\/}/L$.

\begin{lemma}
\label{bijection}
$\Xi: {\cal C\/}/L \to {\cal X\/}/\Lambda$ is a bijection.
\end{lemma}

\startproof
Let $c \in \cal C$.  We have
$$c=(x,y), \hskip 30 pt x=m+1/2, \hskip 10 pt
y=n+1/2, \hskip 10 pt m,n \in \Z.
$$
Looking at the formula
in Equation \ref{map}, we see that
$\Xi(c) \in {\cal X\/}$.

Now ${\cal C\/}/L$ and ${\cal X\/}/\Lambda$
both have $\omega^3$ points.  The idea
is that $[0,\omega^2] \times [0,\omega]$ is
a fundamental domain for $L$ and
$[-1,1]^3$ is a fundamental domain for $\Lambda$.  
Given the equality in the number of points
of the domain and range, we just need to show
that $\Xi$ is an injection from one set into the
other.

Suppose that
$\Xi(c_1)=\Xi(c_2)$.   We write
\begin{equation}
c_2-c_1=(x,y)=\bigg(\frac{m}{\omega},\frac{n}{\omega}\bigg),
\hskip 50 pt m,n \in \Z.
\end{equation}
We have
$$\Psi(c_2)-\Psi(c_2)=
\bigg(\frac{2pm}{\omega^2},
\frac{2pm}{\omega^2}-\frac{pn}{\omega^2},
\frac{2pm}{\omega^2}-\frac{Pn}{\omega^2}\bigg) \hskip 10 pt
{\rm mod\/} \hskip 5 pt \Lambda.
$$
In order for the first coordinate to vanish mod
$\Lambda$, we must have $m \equiv 0$ mod $\omega^2$.
For the second and third coordinates to vanish
mod $\Lambda$, we must have $n \equiv 0$ mod $\omega$.
But then $c_2-c_1 \in L$.
\endproof

So, $\Xi({\cal C\/})$ intersects
$X$ only in the fibers above $T \in \Z_1/\omega$,
and the coordinates of the image are in
$\Z_0/\omega$.

\begin{lemma}
The fiber over $T \in \Z_1/\Omega$ intersects
the walls of the partition in lines of the
form $x=u$ and $y=u$ with $u \in \Z_1/\omega$.
\end{lemma}

\startproof
When $T=-1$, the values of $u$ are $\{-1,1+P,1\}$,
all of which belong to $\Z_1/\omega$.
As the value of $T$ changes by $2/\omega$,
the offsets for the wall-fiber intersections
change by $\pm 2/\omega$.
\endproof

Combining our two results, we see that
$\Xi({\cal C\/})$ never hits a wall
of the partition.  Hence
$(\Xi,X)$ yields a well-defined tiling.

If we place a cube of side length $2/\omega$ around
each point of $\cal X$, then these cubes tile $\R^3$.
Moreover, these cubes intersect each fiber over
$T \in \Z_0/\omega$ in a $\omega \times \omega$ grid
which exactly fills the square $[-1,1]^2$.  The
image $\Xi(\cal C)$ intersects this grid at the
centers of the squares whereas the walls of the
partition intersect the grid in line segments extending
the edges of the squares.

\subsection{The Isomorphism Theorem}

The theorem in this section is a precursor to the
PET Equivalence Theorem. In the next chapter
we will deduce the PET Equivalence Theorem from
the result here.
Referring to the grid description of the plaid model,
we define the {\it grid light set\/} of a square tile $Q$ to
be the set of edges with one light particle in them.
We define the {\it tile light set\/} of a 
square tile to be the
set of edges involved in its connector. 
Here is the main result of the paper.

\begin{theorem}[Isomorphism]
\label{iso}
For any even rational parameter and any
integer unit square, the grid light
set coincides with the tile light set.
\end{theorem}

The Fundamental Theorem for the plaid model is an
immediate consequence.
From the tile description,
the tile light set for any square has $0$ or $2$ members.
Hence, so does the grid light set.
Since the grid and light sets coincide, the
equality of the two systems immediately
implies that the tiling produced by our
classifying pair is coherent.

\subsection{Symmetries of the System}

The tiling produced by our classifing pair is
clearly invariant under the lattice $L$
discussed in the previous section.
Here we show that it is also invariant under
the larger group $\widehat L$ of symmetries
of the plaid model.  As in \S \ref{plaidsymm}
It suffices to prove that reflection in the
origin is a symmetry, and that reflection
in the $x$-axis is a symmetry.

\begin{lemma}
The tiling produced by $(\Xi_P,X)$ is
symmetric with respect to reflection in
the origin.
\end{lemma}

\startproof
We set $\Xi=\Xi_P$ and we drop the
first parameter.  Thus, we think
of $\Xi$ as a map from
$\R^2$ into $\R^3$.
Define
$$\rho(x,y)=(-x,-y),
\hskip 30 pt
\Psi(T,U_1,U_2)=(-T,-U_1,-U_2)$$
It follows directly from the
formulas that
$\Xi(-x,-y)=-\Xi(x,y)$.
At the same time, $\Psi$
preserves our partition of $X$ and permutes
labels according the following
scheme: N and S are swapped and
E and W are swapped.
But that means that $\rho$ maps
the tile centered at $c$ to the tile
centered at $\rho(c)$.
\endproof

\begin{lemma}
\label{refly}
The tiling produced by $(\Xi_P,X)$ is
symmetric with respect to reflection in
the $x$ axis.
\end{lemma}

\startproof
We make the same notational
conventions as in the previous proof.
This time define $\rho(x,y)=(x,-y)$ and
$\Psi(T,U_1,U_2)=(T,U_2,U_1)$.
We compute easily that
$\Xi \circ \rho = \Psi \circ \Xi$,
when these maps are restricted to
points $(x,y)$ with $y$ a half integer.
This set contains $\cal C$.
At the same time, $\Psi$ preserves
our partition (and indeed each fiber)
and permutes the labels by swapping
N and S and doing nothing to E and W.
But that means that $\rho$ maps
the tile centered at $c$ to the tile
centered at $\rho(c)$.
\endproof

\subsection{Some Additional Formulas}
\label{local}

We fix $p/q$ and set $P=2p/(p+q)$ as usual.
For the purposes of calculation, it is 
useful to have formulas for $\Xi_P$ which
take values in the fundamental domain
$\{P\} \times [-1,1]^3$. We use the coordinates
$(T,U_1,U_2)$ discussed above.
Let $[x]_2$ denote the value of $x$ mod $2\Z$ that
lies in $[-2,2)$.

\begin{equation}
\label{loc1}
T(x,y)=[2Px+1]_2.
\end{equation}
\begin{equation}
b(x,y)=\frac{1}{2}PT(x,y).
\end{equation}
\begin{equation}
U_1(x,y)=[PQx+b(x,y)-Py]_2.
\end{equation}
\begin{equation}
U_2(x,y)=[PQx+b(x,y)+Py]_2.
\end{equation}

\begin{lemma}
$\Xi=(P,T,U_1,U_2) \hskip 10 pt
{\rm mod\/}\hskip 5 pt \Lambda$.
\end{lemma}

\startproof
On $\cal C$, the map $T(x,y)$
agrees with $$T'(x,y)=[2Px+2y]_2,$$
because the coordinates of points
in $\cal C$ are odd half-integers.
So, it suffices to check that
$\Xi=(P,T',U_1,U_2)$ mod $\Lambda$.
One can check that right hand
side of this equation
respects the identifications
on $[-1,1]^3$ and gives a
locally affine map into $X_P$.
Next, using the identity $PQx+P^2x=2Px$, one
can check the equation
for an open subset of $\R^2$ in
which $[2Px+1]_2=2Px+1$.  But then
the equality always holds, by
analytic continuation.
\endproof

\subsection{Irrational Limits}
\label{irr}

The classifying pair $(\Xi_P,X)$ makes
sense even when $P$ is irrational.
However, it might not happen that this map
gives a well-defined tiling.  In this
irrational case, 
\begin{equation}
\Xi_P(1/2,1/2)=(-1+P,0,P),
\end{equation}
and this is a point in the boundary of
the partition.

To remedy this situation, we introduce an
{\it offset\/}, namely a vector
$V \in \R^3$, and we define
\begin{equation}
\Xi_{P,V}=\Xi_P+(0,V).
\end{equation}
(The first coordinate, namely $P$, does not change.)
Since $\Lambda_P$ acts on $\R^3$ with
compact quotient, we can take $V$ to lie in a
compact subset of $\R^3$ and still we will
achieve every possible map of this form.
Since the boundary of our partition has measure
$0$, almost every choice of $V$ leads to a
classifying pair that defines a well-defined
tiling.  In this case, we call $V$ a {\it good offset\/}.

\begin{lemma}[Geometric Limit]
Let $V$ be a good offset for $\Xi_P$.
Then the tiling defined by
$(\Xi_P,X)$ is coherent.
\end{lemma}

\startproof
The proof just amounts to taking a limit of the
coherence in the rational case.
Let $\{p_k/q_k\}$ be a sequence of even rational
parameters so that $P_k \to P$.  Here
$P_k=2p_k/(p_k+q_k)$.   Given any vector
$W_k \in \Z^2$, we define
\begin{equation}
\Xi_P^{W_k}(c)=\Xi_{P_k}(c-W_k).
\end{equation}
Since $\Xi_P$ is locally affine, there is
some other vector $V_k \in \R^3$ such that
\begin{equation}
\Xi_{P_k}^{W_k}=\Xi_{P_k,V_k}.
\end{equation}
We can choose the vectors $\{W_k\}$ so that
the vectors $V_k$ converge to our offset $V$.

Suppose that $\tau_1$ and $\tau_2$ are two
adjacent square tiles.  Let $c_1$ and $c_2$
be the corresponding centers.  By definition
$\Xi_{P,V}$ maps both these centers
into the interiors of pieces of our partition.
By continuity, $\Xi_{P_k,W_k}$ also maps these
centers into the interiors of the same pieces,
once $k$ is sufficiently large.  But that
means that the tile types are the same
for $\tau_1$ and $\tau_2$, with respect to
either $P_k$ or $P$.  Since the tiling
is coherent for $(\Xi_{P_k},X)$, the tiles
$\tau_1$ and $\tau_2$ match across their
boundaries.
\endproof

\newpage

\section{The Distribution of Polygons}
\label{loops}
\subsection{Overview}

The first
main goal is to prove the following result.

\begin{theorem}
\label{bigg}
Let $\{p_k/q_k\} \subset (0,1)$ be any sequence of even
rational numbers with an irrational limit.
Let $\{B_k\}$ be any sequence of
associated blocks.
Let $N$ be any fixed integer.
Then the number of distinct plaid polygons
in $B_k$ and the maximum diameter
of a plaid polygon in $B_k$ both exceed $N$
for all $k$ sufficiently large.
\end{theorem}

When we know more
about the sequence, we can get more information.
We say that $\{p_k/q_k\}$ is {\it tuned\/}
if the sequence $\{\tau(p_k/q_k)\}$ also
converges. We call $\lim \tau(p_k/q_k)$ the
{\it tuned limit\/} of the seguence.
By compactness, every sequence has a tuned
subsequence.
We call the tuned sequence {\it (ir)rationally tuned\/}
if the tuned limit is (ir)rational.

We say that the $x$-{\it diameter\/} of a
planar set is the diameter of its projection
to the $x$-axis. This quantity has a dynamical
interpretation; see \S \ref{dyn}.

\begin{theorem}
\label{bigg2}
Let $\{p_k/q_k\} \subset (0,1)$ be any irrationally
tuned sequence. Let $\{B_k\}$ be any sequence of
associated blocks.
Let $N$ be any fixed integer.
Then there is some $\delta>0$ such that
the following property holds once $k$ is sufficiently large:
At least $N$ distinct plaid polygons have $x$-diameter
at least $\delta \omega_k$, and every point of $B_k$
is within $\omega_k/N$ of one of them.
\end{theorem}

\subsection{One Large Polygon}

The result in this section is just a warm-up.
Let $p/q$ be an even rational parameter
and let $\omega=p+q$.

\begin{theorem}
\label{first}
Let $B$ denote the first block.
Then there exists a plaid polygon
in $B$ whose $x$-diameter is at least
$\omega^2/(2q)-1$.  Moreover,
this polygon has bilateral symmetry
with respect to reflection in the
horizontal midline of $B$.
\end{theorem}

\startproof
Let $L$ be the horizontal
line of capacity $2$ and positive sign
which intersects $B$.  Let $z_1=(0,y) \in L$.
By Lemma \ref{anchor}, we
 know that $z_1$ is a light
point of mass $1$.  We compute that
$F_Q(z_2)=1/\omega$, when
$z_2=(\omega^2/2q,y)$.
Hence $z_2$ is another light point on $L$.
Since $L$ has capacity $2$, these are
the only two light points on $L$.
The lattice polygon which crosses the
unit horizontal segment containing $z_1$
must also cross the unit horizontal segment
containing $z_2$ because it has to
intersect $L \cap B$ twice.
This gives the lower bound on the $x$-diameter.

Let
$\Gamma'$ denote the reflection
of $\Gamma$ in the horizontal midline of $B$.
We want to show that $\Gamma'=\Gamma$.
Let $V_1$ and $V_2$ denote the two vertical lines
of $B$ having capacity $2$.  These lines are 
symmetrically placed with respect to the
vertical midline of $B$.  Hence, one of the
two lines, say $V_1$, lies less than $\omega /2$
units away from the $y$-axis.
Since 
$\omega /2<\omega^2/(2q),$
the point $z_2$ is separated from the $y$-axis by $V_1$.
Hence both $\Gamma$ and $\Gamma'$ intersect
$V_1$.  Since there can be at most $1$
plaid polygon which intersects $V_1 \cap B$,
we must have $\Gamma=\Gamma'$.
\endproof

\subsection{The Empty Rectangle Lemma}

Fixing a parameter $p/q$ and a block $B$
and an even integer $K \geq 0$ let
$\Gamma_K$ denote the union of all the
lines of capacity at most $K$ which
intersect $B$.  The complement
$B-\Gamma_K$ consists of $(K+1)^2$ rectangles
arranged in a grid pattern.  We say that
one of these rectangles is {\it empty\/} if
its boundary has no light points on it.
Empty rectangles serve as barriers,
separating the plaid polygons inside
them from the plaid polygons outside them.

\begin{lemma}
For all parameters, all blocks $B$, and all
choices of $K$, at least one of
the rectangles of $B-\Gamma_K$ is empty.
\end{lemma}

\startproof
This is a counting argument.  We will
suppose that there are no empty boxes and
derive a contradiction.  There are a total
of $(K+1)^2$ rectangles.  If some rectangle
$R$ has a light point on it, then it must
have a second light point, because the polygon
$\Gamma$ crossing into $R$ through an edge containing
one of the light points must cross out of
$R$ through another edge.  

The one exceptional
situation is when the light point $z$ lies at
the corner of $R$.  In this
case, one of the edges $E$ of $R$ lies in
a vertical boundary of the block $B$.
Let's consider the case when $E$ lies in
the west boundary of $R$ and $z$ is
the south west corner.  The other cases
are similar.  $\Gamma$
crosses into $R$ through the south edge
of $R$, but then it cannot exit through
$E$ because $E$ lies in the boundary of $B$.
So, even in this exceptional case, there
must be $2$ light points in the boundary of $R$.

If every rectangle has at least $2$ light
points, then there are at least
$(K+1)^2$ light points total.  The idea is
that we double the number of squares but then
observe that we have counted some of the
points twice.

On the other hand, we know that a line of capacity $k$
contains at most $k$ light points.  Since there
are $4$ lines of capacity $k$ for each $k=0,2,...,K$,
this gives a total of
$$8 \sum_{k=1}^{K/2} k = (K+1)^2-1.$$
We have one fewer point than we need.
This is a contradiction.
\endproof

\subsection{No Big Gaps}

\begin{lemma}
\label{gap}
Let $\{p_k/q_k\} \subset (0,1)$ be any sequence of even
rational numbers with an irrational limit.
Then there is some fixed number $R$ with
the following property.
If $p_k \in B_k$ is any point then the disk
of radius $R$ about $p_k$ intersects at least
one plaid polygon in $B_k$ associated to $p_k/q_k$.
\end{lemma}

\startproof
We use the tile description of the plaid model.
Let $P_k=2p_k/(p_k+q_k)$.  Since
$\lim p_k/q_k$ converges, so does
$\lim P_k$.  Let $R=\lim P_k$.
Consider the maps
\begin{equation}
\Xi'_{P_k}(z)=\Xi_{P_k}(z-p_k).
\end{equation}
There are translation vectors $V_k$ so that
\begin{equation}
\Xi'_{P_k}=\Xi_{P_k}+V_k.
\end{equation}
Here $\Xi_{P_k}$ is the map
described in \S \ref{classmap}.
Since the lattice
$\Lambda_{P_k}$ acts on $\R^3$ with
compact quotient, we can take
$\{V_k\}$ to be a bounded sequence of vectors.

By compactness, there is some limit map
\begin{equation}
\Xi'_R=\lim \Xi'_{P_k}.
\end{equation}
The limit map $\Xi'_R$ differs from the
classifying map described in \S \ref{classmap}
by a translation. If this lemma is false, then
$\Xi'_R({\cal C\/})$ is contained entirely 
inside the portion of the partition which
assigns trivial tiles to the squares.  But
on the other hand, from our work in
\S \ref{defined}, we can take a limit and
see that $\Xi'_R({\cal C\/}$ is dense
in the classifying space $X_R$.  This contradicts
the definition of our partition where all
the pieces of the partition have positive 
volume at all parameters.
\endproof

\subsection{Some Congruence Lemmas}

The following lemmas will be helpful in
the next section.

\begin{lemma}
Let $\{p_k/q_k\}$ be a sequence of 
rational numbers with an irrational limit.
Then it is impossible to find uniformly
bounded integers $\alpha_k,\beta_k,\gamma_k$
so that $p_k \alpha_k + q_k \beta_k=\gamma_k$.
\end{lemma}

\startproof
Let $r$ be the limit of our sequence.
Passing to a subsequence, we can assume that
$\alpha,\beta,\gamma$ are independent of $k$.
Dividing by $q_k$ and taking a limit, we
see that
$$r=\lim_{k \to \infty} \frac{\gamma/q_k-\beta}{\alpha}=
-\beta/\alpha.$$
Hence $r$ is rational. This is a contradiction.
\endproof

As usual, let $\omega=p+q$, for any
even rational parameter $p/q$.

\begin{lemma}
\label{growth}
Let $\{p_k/q_k\}$ be a sequence of 
rational numbers with an irrational limit.
Let $a_k$ be such that $2p_ka_k \equiv -1$ mod $\omega_k$.
Then there does not exist a uniformly
bounded integer sequence $\{b_k\}$ so that
$2a_kb_k$ is uniformly bounded mod $\omega_k$.
\end{lemma}

\startproof
Suppose that the sequence
$\{b_k\}$ exists. Then we have
another uniformly bounded integer
sequence $\{c_k\}$ so that
$$2a_kb_k \equiv \hskip 10 pt c_k \hskip 5 pt
{\rm mod\/}\hskip 5 pt \omega_k.$$
Multiplying through by $p_k$ we get
$$-b_k \equiv c_kp_k  \hskip 5 pt
{\rm mod\/}\hskip 5 pt \omega_k.$$
But then
$$c_kp_k = d_k(\omega_k)-b_k,$$
for some uniformly bounded sequence $\{d_k\}$.
Since $\omega_k=p_k+q_k$, we can 
expand this out and we get a contradiction
to the previous result.
\endproof

\subsection{Proof of Theorem \ref{bigg}}

Let $\{p_k/q_k\}$ and $\{B_k\}$ be
as in Theorem \ref{bigg}. 

\begin{lemma}
The diameter of the largest plaid polygon
in $B_k$ tends to $\infty$.
\end{lemma}

\startproof
By symmetry, we can assume that the
block $B$ is one of the fundamental blocks. Let
$H_-$ and $H_+$ denote the two members of
$\cal H$ which have capacity $2$.  Let
$V_-$ and $V_+$ denote the two members of
$\cal V$ which have capacity $2$.  These lines
intersect $B$ in a pattern which has
$4$-fold dihedral symmetry.
We will treat the case when $H_+$ lies to
the left of $H_-$ and $V_+$ lies to
below $V_-$.  The opposite case has
a very similar treatment, and in fact
follows from the first case and symmetry.

The line $H_+$ goes through the point $(0,y)$, where
$2py \equiv 1$ mod $\omega$.   We have
$|y| \to \infty$ by Lemma \ref{growth}.
Hence, the separation between
the lines $H_+$ and $H_-$ tends to
$\infty$ with $k$. By symmetry, the
same goes for the separation between $V_+$ and $V_-$.
In particular, the distance between $H_+$ and
the $x$-axis tends to $\infty$.
It is convenient to set $H=H_+$ and $V=V_+$.

\begin{center}
\resizebox{!}{2in}{\includegraphics{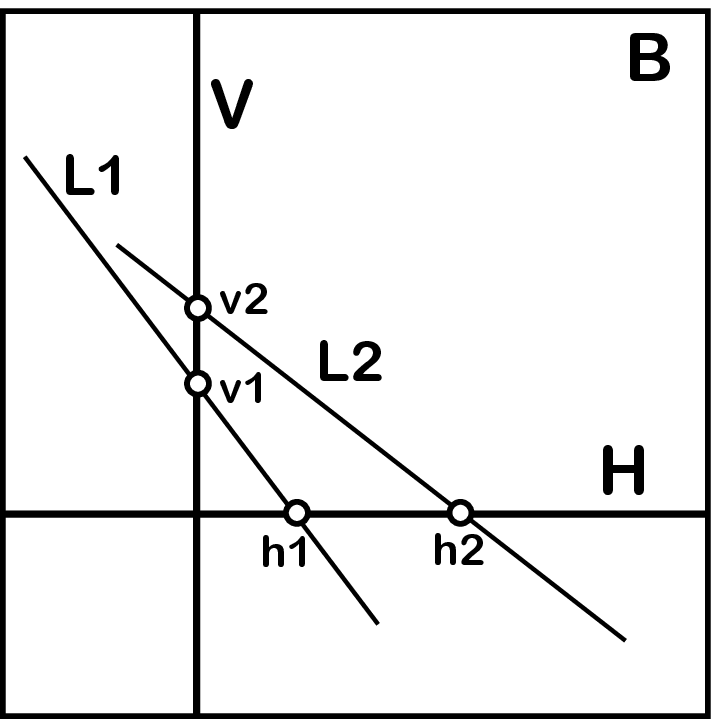}}
\newline
{\bf Figure 4.1:\/}The light points on $H$ and $V$.
\end{center}

Consider the light points $v_1$ and $v_2$
on $V$, as shown in Figure 4.1.
If the distance between these two points
tends to $\infty$ with $k$, then the plaid
polygon through these points has diameter
tending to $\infty$.  So, we just have
to consider the case when the distance
between these points remains uniformly bounded.
Since these points are symmetrically placed
above and below the $x$-axis, both points
must remain within a uniformly bounded
distance of the $x$-axis.

One of the points, say $v_1$, lies on the
$\cal Q$ line $L_1$ of mass $1$.
Since $L_1$ has slope
$-Q<-1$, we know that $L_1$ intersects
$H \cap B$ at one of the light points $h_1$.

In order for us to avoid finishing the
proof of the theorem, the second
light point $h_2$ on $H$ must
remain within a uniformly bounded distance
of $h_1$.  If $h_2$ also lies on
a $\cal Q$ line, then the distance
between the two $\cal Q$ lines of
mass $1$ remains uniformly bounded.
The same reasoning as we gave for
the $\cal H$ lines rules out this
possibility.  Thus, $h_2$ lies on
a $\cal P$ line $L_2$.
Since $L_2$ has slope
$-P \in (-1,0)$, the line
$L_2$ intersects $V \cap B$ in the other
relevant point $v_2$.

Consider the segments of $L_1$ and $L_2$
which have their endpoints on $V$ and $H$.
These two segments tend to $\infty$ in length
and yet each endpoint of $L_1$ is uniformly
close to an endpoint of $L_2$.  At the
same time, the slope of $L_1$ tends
to $-1/P$ and the slope of $L_2$
tends to $-1/Q$.  This situation is
impossible.  This contradiction
finishes the proof.
\endproof

\begin{lemma}
The number of distinct plaid polygons
in $B_k$ tends to $\infty$ as
$k \to \infty$.
\end{lemma}

\startproof
Let $M$ be any fixed positive integer.
The side lengths of the rectangles in the
grid $\Gamma_M$ from the Empty Rectangle
Theorem all have the form $c a_k$ mod $\omega_k$,
where $c$ ranges from a uniformly bounded
set.  By Lemma \ref{growth}, the minimum
side length of a rectangle in $\Gamma_M$
tends to $\infty$ with $k$.

By choosing $k$ sufficiently large, we can
find a collecion $R_{k,1},...,R_{K,N}$ of
rectangles such that each $R_{k,j}$
contains a disk $D_{k,j}$ with the
following two properties.
\begin{itemize}
\item The radius of $D_{k,j}$ tends
to $\infty$ with $k$.
\item $D_{k,j} \cap R_{k,i}=\emptyset$
for all $j< i$.
\end{itemize}
The idea here is that we apply the
Empty Rectangle Lemma at a small scale, then
at a much larger scale, and so on.
In the end, we get $N$ different rectangles
which are either nested or pairwise disjoint,
and in either case we can choose the scales
to get the big disks.

It follows from Lemma \ref{gap} that each $D_{k,j}$
of these disks intersects a plaid polygon $\Gamma_{k,j}$
once $k$ is sufficiently large.
By construction, the plaid polygons
$\Gamma_{k,1},...,\Gamma_{k,N}$ are
pairwise distinct because
they are separated by the rectangles.
\endproof

\subsection{Rescaled Limits}

Now we turn to the case when $\{p_k/q_k\}$ is irrationally
tuned. Let $\tau \in (0,1)-\Q$ be the tuned limit.
 There is a unique homothety $T_k$ such that
$T_k(B_k)=[0,1]^2$.   In this case, for each value of
$N$, the rescaled grids
\begin{equation}
T_k(\Gamma_{k,N})
\end{equation}
converge (in the Hausdorff topology, say) to a grid
$G_M$ which has the following description.
\begin{itemize}
\item For each $k=0,...,M$, the grid $G_M$
contains the horizontal lines
$y=k\tau$ and $y=1-k\tau$.  These quantities are
taken mod $1$.  These are rescaled limits
of the horizontal lines of capacity $k$.
\item For each $k=0,...,M$, the grid $G_M$
contains the vertical lines
$x=k\tau$ and $x=1-k\tau$. These quantities are
taken mod $1$.  These are rescaled limits
of the vertical lines of capacity $k$.
\item $G_M$ has no other lines.
\end{itemize} 

Now we consider the $\cal P$ lines and
the $\cal Q$ lines. The situation
here is more subtle.
\begin{equation}
P=\lim P_k, \hskip 30 pt
P_k=\frac{2p_k}{\omega_k},
\hskip 30 pt Q=2-P.
\end{equation}
To understand the subtlety, let's first consider
the case when $B_k$ is the first block for all $k$.

\begin{lemma}
Let $M$ be a positive odd integer.
The rescaled limit
$T_k(\Omega_{k,M})$ exists.  It consists of
the lines of slope $-P$ and $-Q$ which
have $y$-intercept
$\pm \mu \tau + \lambda$ for
$\mu \in \{1,3,5,...,M\}$ and
$\lambda \in \Z$.
\end{lemma}

\startproof
A calculation very much like the one in
the proof of Lemma \ref{anchor} shows that,
for $\mu$ odd and $\lambda \in \Z$,
\begin{equation}
F_{P_k}(0,\pm \mu a_k+\lambda \omega)=
F{P_k}(0,\pm \mu)=\mp \mu/\omega_k.
\end{equation}
Hence the $\cal P$ and $\cal Q$ lines through
$(0,\pm \mu)$ have mass $\mu$.
Once $k$ is sufficiently large,
all these lines belong to $\Omega_{k,M}$ for
$\mu=1,3,5,..,M$.  Moreover, no other
points of the form $(0,y)$ satisfy
$2py \equiv \pm \mu+\omega$ mod $2\omega$.
\endproof

Now we consider the case when $\{B_k\}$ is an
arbitrary sequence of blocks. Let
$B_k^0$ denote the first block associated to
$p_k/q_k$ and let
$\Omega_{k,M}^0$ denote the corresponding
set of lines.
The set $T_k(\Omega_{k,M})$ differs from
the set $T_k(\Omega_{k,M})$ by a vertical
translation.  We can take this translation
to have length at most $1$ because our
sets are both invariant under vertical translation
by $1$.  So, by compactness, we can pass to a
subsequence and assume that
$T_k(\Omega_{k,M})$ really does converge.  The
limit is the set of lines of slope $-P$ and $-Q$ 
having $y$-intercept $\pm \mu + \lambda + \xi$,
where $\mu \in \{1,3,5,...,M\}$ and
$\lambda \in \Z$ and $\xi \in (0,1)$ is the
translation factor.   We call $\xi$
the {\it offset\/} of the limit.

So, if we pass to a subsequence, then the sets
$\{T_k(\Gamma_{k,M})\}$ and
$\{T_k(\Omega_{k,M})\}$ converge
to $\Gamma_M$ and $\Omega_M$ respectively.
(In the final section of this chapter, we
will describe a more canonical way to take
rescaled limits.)
We assign a capacity to the lines
in $\Gamma_M$ in the obvious way:  If some
line is the limit of lines of capacity $c$,
it gets capacity $c$.  Likewise, we assign
a mass to the lines in $\Omega_M$.
This allows us to assign a set of light
$\Sigma_M$ on the lines of $\Gamma_M$.

By construction, the sets
$T_k(\Sigma_{k,M})$ {\it practically converge\/}
to the set $\Sigma_M$.  There is one case we
have to worry about.  If $\Sigma_M$ contains a
point in the corner of $[0,1]^2$ then it
might not arise as the limit of points in
$T_k(\Sigma_{k,M})$.  This situation would not
happen if $B_k$ is always the first block,
but it could happen in general.  But we
can say that every point of $\Sigma_M \cap (0,1)^2$
is the limit of points in $T_k(\Sigma_{k,M})$.
Moreover, if a light point in $\Sigma_M \cap (0,1)^2$
is at least $\delta$ from every other point
of $\Sigma_M$, then the corresponding
point of $\Sigma_{k,M}$ is separated by
$\delta\omega_k$ from every other light
point in $\Sigma_{k,M}$.

The really interesting thing is that it seems that we
can make sense of rescaled limits of individual
polygons.  This is beyond the scope of this paper,
but we will discuss it in the last section.

\subsection{The Filling Property}

Now we recall a familiar fact about circle
rotations.  Let $\theta \in (0,1)$ be
irrational.  Consider the map
\begin{equation}
T(x)=x+\theta \hskip 5 pt {\rm mod\/} \hskip 5 pt 1.
\end{equation}
For any $\epsilon>0$ there is some
$M$ such that the first $M$ points of the orbit
$\{T^j(x)\}$ is $\epsilon$-dense.  The value of
$M$ only depends on $\theta$ and
$\epsilon$ and not on the
starting point $x$.  The way that
$M$ depends on $\theta$ and $\epsilon$
is subtle; it has to do with the
continued fraction expansion of $\theta$.
However, we do not care about this
subtlety.

Here is a consequence of the filling property.
We keep the notation from the previous section. 

\begin{lemma}
Let $D$ be a disk of radius $\epsilon$ in $[0,1]^2$.
Then there is a constant $M$ and some
horizontal line $L$ of $\Gamma_M$ so that
$L \cap D \cap \Sigma_M$ contains at least
$2$ points.
\end{lemma}

\startproof
We say that a line $L$ {\it frankly intersects\/}
a disk $D$ if $L \cap D$ contains a point which is within
${\rm radius\/}(D)/4$ from the center of $D$.
If a horizontal line and a line of slope $-Q \in (-1,-2)$
both frankly intersect $D$, then their intersection
is contained in $D$.

By the filling property, there is some
$M'$ such that at least $2$ lines of
$\Omega_{M'}$ frankly intersect $D$. 
Call these lines $Q_1$ and $Q_2$.
We can take these lines to be of type $\cal Q$
and of positive sign.
Also by the filling property, there is some
$M>M'$ so that at least one horizontal
line $L$ having positive sign and capacity
in $(M',M)$ frankly intersects $D$.

We get out two points of $L \cap D \cap \Sigma_M$
by intersecting $L$ with the two lines $Q_1$ and $Q_2$.
\endproof

We say a bit more about the points produced by the
previous lemma.

\begin{lemma}
Suppose that $z_1$ and $z_2$ are two
points of $\Sigma_M \cap D \cap \Gamma_M$
that lie on the same horizontal line.
Then at least one of the two points is
distinct from every other point of
$\Sigma_M \cap L$.
\end{lemma}

\startproof
Let $L$ be the line $y=\xi_1$.
Let $\xi_2$ be the offset of our limit.
Our points
have the form
\begin{equation}
z_j=\bigg(\frac{\mu_j \tau + \xi_3}{Q}\bigg), \hskip 30 pt
|\mu| \leq M, \hskip 30 pt
\xi_3=\xi_1+\xi_2.
\end{equation}
Given the irrationality of $\tau$, these two points
are distinct from each other, and also distinct
from every other point of type Q on $L$.

Suppose then that $z_1$ and $z_2$ are both points of
type P as well.   Then we have
\begin{equation}
|z_1-z_2|=\frac{c_1}{P}=\frac{c_2}{Q}, \hskip 30 pt
c_1,c_2 \in \Z.
\end{equation}
This contradicts the fact that $P/Q$ is irrational.
\endproof

\subsection{Proof of Theorem \ref{bigg2}}

We will work with the rescaled limit and then,
at the end, interpret what our result says.

The Empty Rectangle Lemma applies to the
grid $\Gamma_M$  By varying $M$ and applying
the Empty Rectangle Lemma $N$ times, we
can find $N$ rectangles in $[0,1]^2$,
say $R_1,...,R_N$ such that 
each $R_j$ contains a disk $D$ which
is disjoint from $R_i$ for $j<i$.
The filling property lets us take 
$N$ as large as we like.  

By the work in the previous section, each
disk $D_j$ will contain a light point $z_j$.
which is distinct from all other light
points on the same horizontal line.
Let $2\delta$ be the minimum
separation between $z_j$ and any other
light point.  The minimum is taken
over all $j=1,...,N$.  
By making $\delta$ smaller, if necessary,
we can arrange that every disk of radius
$1/2N$ contains the kind of light point
which is separated from its horizontal
neighbors by at least $2\delta$.

Now let us see what this says about
the picture in the block $B_k$.  
Once $k$ is sufficiently large, we
can find rectangles
$R_{k,1},...,R_{k,N}$ which contain
disks $D_{k,1},...,D_{k,N}$ having
the following propertyes.
\begin{itemize}
\item $D_{k,j}$ is disjoint from $R_{k,i}$ if $j<i$.
\item $D_{k,j}$ contains a light point $z_{k,j}$
which is separated from all other light points
on the same horizontal line by at least
$\delta \omega_k$.
\end{itemize}
Let $\Gamma_{k,j}$ denote the plaid polygon
that intersects the horizontal unit segment
containing $z_j$.  Since
$\Gamma_{k,j}$ is a closed loop, it must
intersect the horizontal line containing
$z_{k,j}$ in a second light point.
Hence $\Gamma_{k,j}$ has $x$-diameter
at least $\delta \omega_k$.
By construction, the polygons
$\Gamma_{k,j}$ and $\Gamma_{k,i}$ lie
in different components of
$B_k-\partial R_{k,j}$ for $j<i$.
Hence, these polygons are all distinct.

At the same time, choose any tile center
$c \in B_k$.  Let $\Delta$ be the disk
of radius $\omega_k/N$ about $c$.  Then
$T_k(\Delta)$ contains a disk of radius
$1/2N$ for $k$ large.  Hence
$T_k(\Delta)$ contains a light
point that is separated from its
horizontal neighbors by at least $2\delta$.
But then the inverse image of this point is a
light point in $\Delta$ that is
is separated from all the other light
points on the same horizontal line by
at least $\delta \omega_k$.  Hence,
$c$ is within $\omega_k/N$ of some
plaid polygon having diameter at least
$\delta \omega_n$.

This completes the proof of
Theorem \ref{bigg2}.

\newpage

\section{The Action on Particles}
\label{PARTICLE}

\subsection{Overview}

Let $\cal C$ denote the set of centers of
square tiles. Let $\Xi$ be the classifying
map from \S \ref{classmap}. To prove the
Isomorphism Theorem we need to study
the image $\Xi({\cal C\/})$.
One could say that our proof of
the Isomorphism Theorem has a product
structure to it.  In this chapter, we underStand
what the classifying map does to those
subsets of $\cal C$ corresponding to
the particles discussed in \S \ref{particle}.
We won't analyze where these subsets
get mapped in, but we will be able to
describe their geometry.  That is, we
will understand their images up to translation.

In the next chapter, we understand where
precisely the classifying map takes
certain elements of $\cal C$ which
correspond to specially chosen instances
of the particles.  This information allows
us to anchor the geometric picture obtained
in this chapter, so to speak.  In \S 7 we
will put the two pieces of information
together and finish the proof of the
Isomorphism Theorem.
In our analysis,
we will work with the west and south edges.
It turns out that this is all we need.

We fix an even rational $p/q \in (0,1)$.
Let $\omega=p+q$.
  Let $z_0,...,z_{n-1}$ denote
the portion of a particle 
in which $z_j$ does
not change type.  There are three cases.
\begin{itemize}
\item  If the particle is vertical,
then $n=\omega$.
\item  If the particle is
horizontal and the points have type P
then $z_0$ is contained in the left edge
of the boundary and $n=2p$.

\item  If the particle is
horizontal and the points have type Q
then $z_0$ is contained in the right edge
of the boundary and
$n=2q$.
\end{itemize}

Let $c_0,...,c_{n-1}$ denote the set of centers of
unit squares in the block such that
$z_i$ belongs to the west (respectively south)
edge of the square centered at $c_i$ provided
that the particle is vertical (respectively horizontal).
We are interested in the image
\begin{equation}
Z =
\bigcup_{i=0}^{n-1} \Xi(c_j) \subset X_P.
\end{equation}
The goal of this chapter is to understand the
set $Z$.

\subsection{The Vertical Case}

For our analysis, we will drop off the first coordinate
of $\{P\} \times X_P$ and work in $\R^3$.
Recall that $X_P$ has
coordinates $(T,U_1,U_2)$, with 
$(U_1,U_2)$ being the fiber.

\begin{lemma}
Suppose that $Z$ corresponds to a vertical
particle of type P.  Then
$Z$ is contained in a single fiber and
lies in a line $U_1={\rm const.\/}$
\end{lemma}

\startproof
Let $c=c_j$ and $c'=c_{j+1}$.  
We will suppose that the
particle is of type P.
This means that 
\begin{equation}
c'=c+(a\omega,1), \hskip 30 pt
2ap \equiv -1 \hskip 10 pt \hskip 5 pt {\rm mod\/} \hskip 5 pt \omega.
\end{equation}
According to Equation \ref{map} we have
\begin{equation}
\Psi(c')-\Psi(c)=(2Pa\omega,2Pa\omega-P,2Pa\omega+P) \hskip 10 pt
{\rm mod\/} \hskip 5 pt \Lambda_P.
\end{equation}
Note that $Pa\omega \in \Z$.  Subtracting off
$Pa\omega$ times $(2,2+P,2+P)$, we get
\begin{equation}
\label{vertP}.
\Psi(c')-\Psi(c)=(0,P(-1-Pa\omega),P(1-Pa\omega)) \hskip 10 pt
{\rm mod\/} \hskip 5 pt \Lambda_P.
\end{equation}
Note that
$-1-Pa\omega=K\omega$ for some $K \in \Z$.
Hence
$K\omega P=2Kp \in 2\Z.$
Hence the $U_1$ coordinate in
Equation \ref{vertP} vanishes.
\endproof

The same argument, or an appeal to symmetry,
yields the following.

\begin{lemma}
\label{QVERT}
Suppose that $Z$ corresponds to a vertical
particle of type Q.  Then
$Z$ is contained in a single fiber and
lies in a line $U_2={\rm const.\/}$
\end{lemma}

\subsection{A Revealing Picture}

Before we deal with the horizontal particles,
we draw a picture of $Z$ in the
$(T,U_j)$ plane.  Up to vertical translation,
the images all look the same.  They
are all contained in parallelograms
whose sides have slope $0$ and $1$.
These parallelograms look the same for
both $j=1$ and $j=2$, up to translations.
The left hand side pertains to the
type P portion of the particle
and the right hand side pertains to the
type Q portion. The middle
picture shows how they two pieces fit
together, and also labels some distances
in the picture.

\begin{center}
\resizebox{!}{1.6in}{\includegraphics{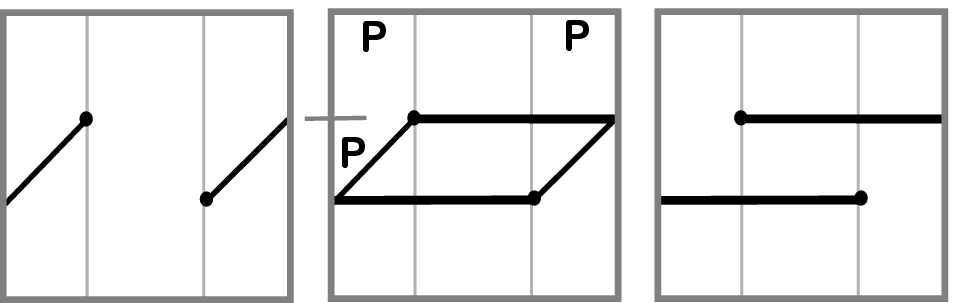}}
\newline
{\bf Figure 5.1:\/} A typical set $Z$ drawn
in the $(T,U_j)$ plane.
\end{center}

It is worth pointing out that the parallelogram
in the middle figure is not embedded.
It has a double point which appears when
the left edge of the picture is identified
to the right edge by the map
$(-1,y) \to (1,y+P)$.  This double point
is the image of the midpoint of one horizontal
edge that is right in the middle of a block.
All the particles pass twice through
a double point like this.

\subsection{The Horizontal Case: Type P}

\begin{lemma}
\label{typeP}
Suppose that $Z$ corresponds to a horizontal
particle of type P.  Then $Z$ is contained in
a single line segment $S$ parallel to
$(1,1,1)$ which does not intersect and of the
fibers over $T \in (-1+P,P-1)$.
\end{lemma}

\startproof
To simplify the argument, we claim that the
result is true for one particle 
$z_0,...,z_{2p-1}$ of this type if and only if
it is true for another particle
$z_0',...,z_{2p-1}'$ of this type.
Let $c_0',...,c_{2p-1}'$ and $Z'$ be the
objects associated to this other
particle.  Since both $z_1$ and $z'$ belong
to the left edge of some block, we have
\begin{equation}
z_j'=z_j+(k_1(p+1),k_2), \hskip 10 pt
c_j'=c_j+(k_1(p+1),k_2), \hskip 30 pt
\forall j.
\end{equation}
\begin{equation}
\Xi(c_j')-\Xi(c_j) = (2Pk_1\omega,\alpha_1,\alpha_2=
 (0,\beta_1,\beta_2) \hskip 20 pt
{\rm mod\/}\hskip 10 pt \Lambda.
\end{equation}
This works because $2Pk_1\omega$ is an even integer.
Here $\alpha_1,\alpha_2,\beta_1,\beta_2$ are
independent of $j$.  Thus $Z'$ is obtained from
$Z$ by applying some fiber-preserving transformation.
This establishes our claim.

In view of our claim, it suffices to work with the particle
\begin{equation}
z_k=(ka\omega+k/P,0), \hskip 30 pt k=0,...,2p-1.
\end{equation}
As usual, $2ap \equiv -1$ mod $P$.

To simplify our calculation, we work with the
points $\widehat c_0,...,\widehat c_{2p-1}$,
where $\widehat c_j$ lies on the midpoint of
the horizontal
segment containing $z_j$. We have
$\widehat c_0=c_0-(0,1/2)$.
Furthermore, we work with the map
$\widehat \Xi=\Xi-(1,0,0)$.  We have
\begin{equation}
\widehat \Xi(x,0)=(2Px,2Px,2Px).
\end{equation}
Let $\widehat Z$ be the modified version of $Z$.
Using the formulas in \S \ref{local}, we see that
the statement of this lemma is equivalent to the
statement that $\widehat Z$ is contained in a line
segment parallel to $(1,1,1)$ and 
is contained in the union of fibers lying over
$[-P,P]$.

Note that $2P(k/P)$ is always an even integer.
Hence $\widehat \Xi(z_j)=(0,0,0)$.  The
points $\widehat c_j$ and $z_j$ are on the same
unit segment and differ by at most $1/2$ units.
Hence $\widehat \Xi(\widehat c_j)$ lies
on the segment joining $(-P,-P,-P)$ to
$(P,P,P)$.  
\endproof

\subsection{A Technical Result}

Recall that $[x]_2$ is the
value of $x$ mod $2\Z$ that lies in $[-2,2)$.
Recall also that $P+Q=2$.  
Here we prove a technical result.

\begin{lemma}
\label{tech0}
Let $m,k \in \Z$ and let $c=m+1/2$.
If $k/Q \in [m+1/Q,m+1]$ then
$[2Pc+1]_2 \in [-1+P,1-P]$.
\end{lemma}

\startproof
The conditions imply that
$2k \in [2Qm+2,2Qm+2Q]$.  Hence, there
is another even integer
integer 
$$e \in [-2Qm-(2Q-2),-2Qm]=[-2Qm-(2-2P),-2Qm].$$
We have
$$2Pc+1=2Pm+P+1 \equiv -2Qm-(1-P) \hskip 5pt
{\rm mod\/} \hskip 5 pt 2\Z.$$
But then
$$[2Pc+1]_2=-2Qm-(1-P) - e \in [-1+P,1-P].$$
This completes the proof.
\endproof

\subsection{The Horizontal Case: Type Q}

\begin{lemma}
\label{typeQ}
Suppose that $Z$ corresponds to a horizontal
particle of type Q.  Then $Z$ is contained in
a single line segment $S$, parallel
to $(1,0,0)$, with the following
structure.
\begin{itemize}
\item $S$ intersects the fibers over
$[-1,-1+P) \cup [P_1,1]$ once.
\item $S$ intersects the fibers
over $[-1+P,1-P]$ twice.
\end{itemize}
\end{lemma}

\startproof
By the same kind of argument given in the
previous result, it suffices to work with the
dark particle $z_0,...,z_{2q-1}$, where
\begin{equation}
z_k=(ka\omega-k/Q,0).
\end{equation}
Likewise, we work with the
points $\widehat c_0,...,\widehat c_{2q-1}$.
This time we will work with the map $\Xi$.
There are $2$ cases:
\newline
\newline
{\bf Case 1:\/}
Suppose that $c=m+1/2$ and
$z_k \in [m+1/Q,m+1]$ for some integer $m$.
By Lemma \ref{tech0}, we can say that
$\Xi(\widehat c_k)$
lies in a fiber over $[-1+P,1-P]$.
In this case there are $2$ points of type Q in
the interval containing $\widehat c_k$ and so
\begin{equation}
\widehat c_{k+1}=\widehat c_k + (a\omega,0).
\end{equation}
Since $2Pa\omega=4ap$, we have
$$\Xi(\widehat c_{k+1})-\Xi(\widehat c_k)=
(4ap,4ap,4ap)=(0,-2Pap,-2Pap) \hskip 10 pt {\rm mod\/}\hskip 5pt \Lambda_P.$$
But $-2Pap \equiv P$ mod $2\Z$.  Hence
\begin{equation}
\label{case1Q}
\Xi(\widehat c_{k+1})-\Xi(\widehat c_k)=(0,P,P) \hskip 10 pt
{\rm mod\/}\hskip 5 pt \Lambda_P.
\end{equation}

\noindent
{\bf Case 2:\/}
Now suppose that $c=m+1/2$ and
$z_1 \in (m,m+1/Q)$.
(The case $z_k=m$ does not occur; it
would be the point $z_{2q}$.)
In this case, we have
\begin{equation}
\widehat c_{k+1}=\widehat c_k+(a\omega-1,0).
\end{equation}
It follows from Equation \ref{case1Q} that
\begin{equation}
\Xi(\widehat c_{k+1})-\Xi(\widehat c_k)=(-2P,P,P)=(2-2P,0,0) \hskip 10 pt
{\rm mod\/}\hskip 5 pt \Lambda_P.
\end{equation}

Let $Z_k=\Xi(\widehat c_j)$.
We have $Z_0=(1,P/2,P/2)$.
In Figure 5.1 this point would be all the
way on the right.
Let $L$ be the line segment from
the lemma which contains $Z_0$.
If we forget the identifications on
$[-1,1]^3$, then $L$ becomes $2$ line
segments, the {\it upper half\/} contains
$Z_0$ and the {\it lower half\/} does not.
At the same time, we divide
$[-1,1]^3$ into $3$ {\it zones\/},
depending on where $T$ falls in 
the partition $[-1,-1+P,1-P,1]$.
The vertical line segments in Figure 5.1 demarkate
these zones.
We say that a point is {\it left lower\/} if it
lies in the left zone and on the lower half of $L$.
We give similar names to the other possibilities.
The two cases above tell us the following.
\begin{itemize}
\item If $Z_k$ is left lower, then
$Z_{k+1}$ is left lower or middle lower.
\item If $Z_k$ is middle lower then
$Z_{k+1}$ is middle upper.
\item If $Z_k$ is middle upper then
$Z_{k+1}$ is right upper or left lower.
\item If $Z_k$ is right upper then
$Z_{k+1}$ is left lower or middle lower.
\end{itemize}
In the last two cases, we use the action of
$\Lambda$ in case the addition of $2-P$ to
the $T$ coordinate causes the value to
increase beyond $1$.
Since $Z_0 \in L$, we have $Z_k \in L$ for all $L$.
\endproof

For later use, we want to extract something more
out of the proof of Lemma \ref{typeQ}.
The following result is a consequence of our
argument given in the proof of Lemma \ref{typeQ}.

\begin{lemma}
\label{typeQQ}
Suppose that $z_1$ and $z_2$ are two intersection
points contained in the same south edge of some unit
square centered at the point $c$.  Suppose that $z_1$
lies to the left of $z_2$.  Then there are
other instances $z_1'$ and $z_2'$ of $z_1$ and $z_2$
respectively so that
\begin{itemize}
\item $\Xi(c_1)$ lies in the fiber over $T=1$ and
$\Xi(c_1)$ and $\Xi(c)$ are joined by a segment
in $[-1,1]^3$ parallel to $(1,0,0)$.

\item $\Xi(c_2)$ lies in the fiber over $T=-1$ and
$\Xi(c_2)$ and $\Xi(c)$ are joined by a segment
in $[-1,1]^3$ parallel to $(1,0,0)$.
\end{itemize}
Here $c_j \in \cal C$ is a point such that the
south edge of the square containing $z_j'$ is
contained in $c_j$.
\end{lemma}

\newpage

\section{The Images of Symmetric Points}
\label{instance}

\subsection{Overview}

We fix an even rational $p/q$ as usual.
Let $\omega=p+q$.   Let $z_1,...,z_n$ be
a particle and let $c_1,...,c_n$ be
the edges corresponding to these
points.  If $\{z_i\}$ is a vertical
(respectively horizontal) particle,
then $c_i$ is the center of the unit square
whose west (respectively south) 
edge contains $z_i$.

In the vertical case, we call $z_i$ a {\it symmetric instance\/}
of the particle if $c_i$ is centered on the horizontal
midline of the block.  That is 
\begin{equation}
\label{symmV}
c_i=(x+1/2,k\omega+\omega/2), \hskip 30 pt
k,x \in \Z.
\end{equation}
In the horizontal case, we call $z_i$
a symmetric instance if $c_i$ is centered on the
vertical midline of the block.  That is, 
\begin{equation}
\label{symmH}
c_i=(k\omega+\omega/2,y+1/2), \hskip 30 pt
k,y \in \Z.
\end{equation}
Every vertical particle has one symmetric instance,
and every horizontal particle has two symmetric
instances - one as type P and one as type Q.  In
the horizontal case, the symmetric instances of
the particle correspond to the double point in Figure 5.1.

Let ${\cal V\/}_k$ denote the subset of
points $c$ in Equation \ref{symmV} such that
the west edge of the
square centered at $c$ has $k$ light
points of type P in them.  Likewise,
let  ${\cal H\/}_k$ denote the subset of
points $c$ in Equation \ref{symmH} such that
the soth edge of the 
square centered at $c$ has $k$ light
points of type P in them.  These
sets are empty unless $k=0,1$.

In this chapter, we will 
$\Xi({\cal V\/}_k)$ and 
$\Xi({\cal H\/}_k)$ 
for $k=0,1$.  (By symmetry, these sets are
the same if we use type Q particles instead
of type P particles.)
Once we have this information,
we will combine it with what we
did in the last chapter to finish the
proof of the Isomorphism Theorem.

For our analysis, it is awkward to work
directly with the centers of the squares.  It
is better to work with the centers of the
edges which contain the particles.  We define
\begin{equation}
{\cal V\/}'_k={\cal V\/}_k-(1/2,0), \hskip 40 pt
{\cal H\/}'_k= {\cal H\/}_k-(0,1/2).
\end{equation}
In each case, we will analyze what
$\Xi$ does to these sets, and then
we will 
make the tedious but routine translation back
to the original sets of interest to us.

\subsection{Good Pairs}

Suppose that $a_1,a_2$ are real numbers.
We say that $(a_1,a_2)$ is a {\it good pair\/} if
none of the numbers
$a_1,a_2,a_1-a_2$ is an integer and
\begin{equation}
[a_1]_2\ [a_2]_2>0, \hskip 30 pt
||a_1]_2|<|[a_2]_2|.
\end{equation}
These are exactly the conditions we used
in the grid description of the plaid model.
Figure 6.1 shows the set of good pairs $(a_1,a_2)$ in
the square $[-1,1]^2$.

\begin{center}
\resizebox{!}{2.5in}{\includegraphics{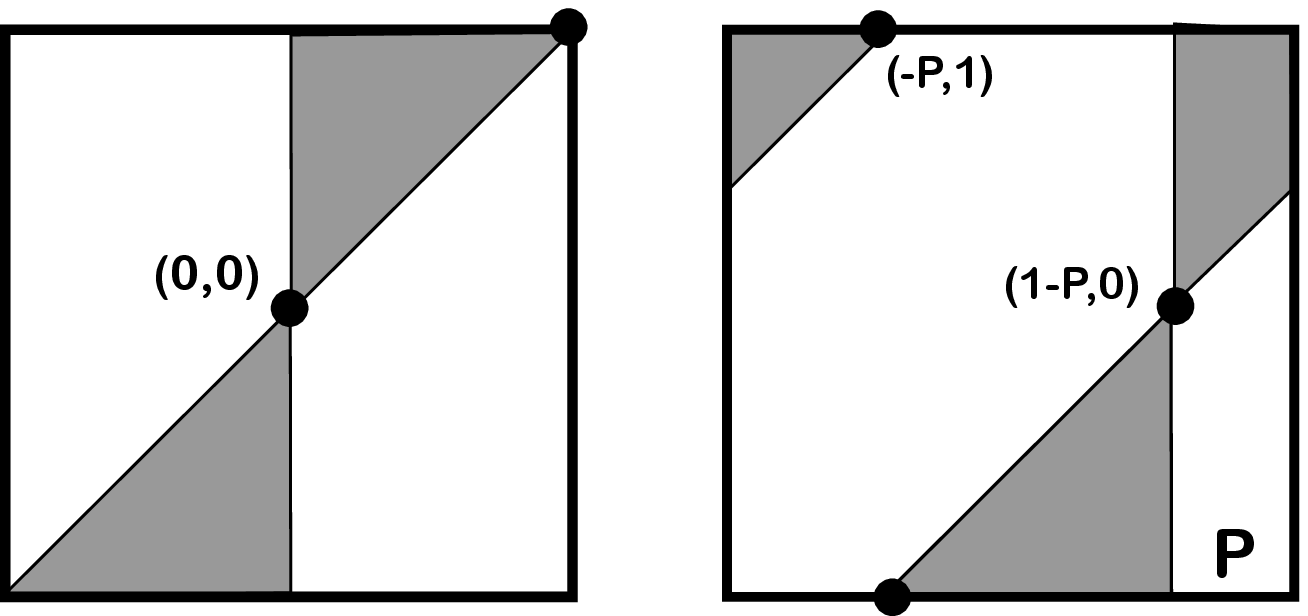}}
\newline
{\bf Figure 6.2:\/} The set of good pairs
\end{center}

Now we prove a technical lemma about good pairs
which will be helpful in the horizontal case.
Define
\begin{equation}
b_1=a_1-a_2, \hskip 30 pt
b_2=a_1-1.
\end{equation}

\begin{lemma}
\label{SWAP}
$(a_1,a_2)$ is a good pair
if and only if $(b_1,b_2)$ is a
good pair.
\end{lemma}

\startproof
Since the affine transformation
$T(x,y)=(x-y,x-1)$ preserves $\Z^2$,
the pair $(a_1,a_2)$ satisfies the
non-integrality condition if and only
if the pair $(b_1,b_2)$ does.
If this lemma is true for the inputs
$(a_1,a_2)$ it is also true for
the inputs $(a_1+2k_1,a_2+2k_2)$
for any integers $k_1,k_2$.  For
this reason, it suffices to consider
the case when $a_1,a_2 \in (-1,1)$.
From here, an easy case-by-case
analysis finishes the proof.
For instance, if $0<a_1<a_2$ then
$0>b_1>b_2$.  The other cases
are similar.
\endproof

\subsection{The Vertical Case}

We drop the first coordinate and work
in $[-1,1]^3$ with the $(T,U_1,U_2)$
coordinates.
We define the {\it diagonal\/} $\Delta X$
of $X$ to be the set where $U_1=U_2$.

\begin{lemma}
\label{QVERT2}
Let $G' \subset \Delta X$ be the
set such that $(T,U_k)$ is a good pair.
Then $\Xi({\cal V\/}'_0) \subset \Delta X-G'$ and
$\Xi({\cal V\/}'_1) \subset G'$ and
\end{lemma}

\startproof
Let $c\in {\cal V\/}'_k$.  
It follows from 
Lemma \ref{refly} (or direct calculation) that
$\Xi(c) \in \Delta X$.
The corresponding intersection point $z=(x,y)$ lies
on a line of slope $-P$ which
has an integer $y$-intercept. Hence
$$
y=\omega/2 - [Px+\omega/2]_1.
$$
Here $[t]_1$ denotes the number in $[-1/2,1/2)$ that differs
from $t$ by an integer. Since
$[t]_1=[2t]_2/2$, we get
\begin{equation}
y=\omega/2-\frac{1}{2}[2Px+\omega]_2=^*
\omega/2-\frac{1}{2}[2Px+1]_2=\omega/2-T(c)/2.
\end{equation}
The starred equality uses the fact that $\omega-1 \in 2 \Z$.
Define
\begin{equation}
z'=(x,y)+T(c)/2.
\end{equation}
Reflection in the horizontal midline of the block
swaps $z$ and $z'$.  By symmetry, $z$ is a light
point of type P if and only if $z'$ is a light
point of type Q.  We will work with $z'$.
We compute
\begin{equation}
\label{great1}
F_V(z')=[2Px]_2=[T(c)+1]_2.
\end{equation}
Next, we compute
$$F_Q(z')=[Py+PQx+1]_2=[P\omega/2+PT(c)/2+PQx+1]_2.$$
Looking at the equations from \S \ref{local},
we recognize this last expression as $[U_2(c)+1]_2$. Since
$U_2(c)=U_1(c)$, we have
\begin{equation}
\label{great2}
F_Q(z')=[U_1(c)+1]_2.
\end{equation}
From Equations \ref{great1} and \ref{great2} we see
that $F_V(z')$ and $F_Q(z')$ have the same
sign if and only if $T(c)$ and $U_1(c)$ have
the same sign and
$|F_Q(z')|<|F_V(z')|$ iff
$|T(c)|<|U(c)|$.
Hence $\Xi(c) \in \Delta X-G'$ if $k=0$ and
$\Xi(c) \subset G'$ if $k=1$
\endproof

\subsection{The Vertical Case Translated}
\label{vtrans}

Now we translate the picture so that we are
working with $\Xi({\cal V\/}_k)$ rather than
$\Xi({\cal V'\/}_k)$.  Let
\begin{equation}
\rho(x,y)=(x+1/2,y).
\end{equation}
We have
\begin{equation}
\rho({\cal V'\/}_k)={\cal V\/}_k.
\end{equation}

A direct calculation shows that
\begin{equation}
\Psi \circ \Xi = \Xi \circ \rho, \hskip 30 pt
\Psi(x_1,x_2,x_3)=(x_1,x_2,x_3)+(P,P,P).
\end{equation}
Hence
\begin{equation}
\Xi({\cal V\/}_k)=
\Xi({\cal V'\/}_k)+(P,P,P) \hskip 30 pt
{\rm mod\/}\hskip 5 pt \Lambda_P.
\end{equation}
Note that $\Psi(\Delta X)=\Delta X$. 
We draw $G'$ and $G=\Psi(G)$ side by side
in Figure 6.2.

\begin{center}
\resizebox{!}{2.2in}{\includegraphics{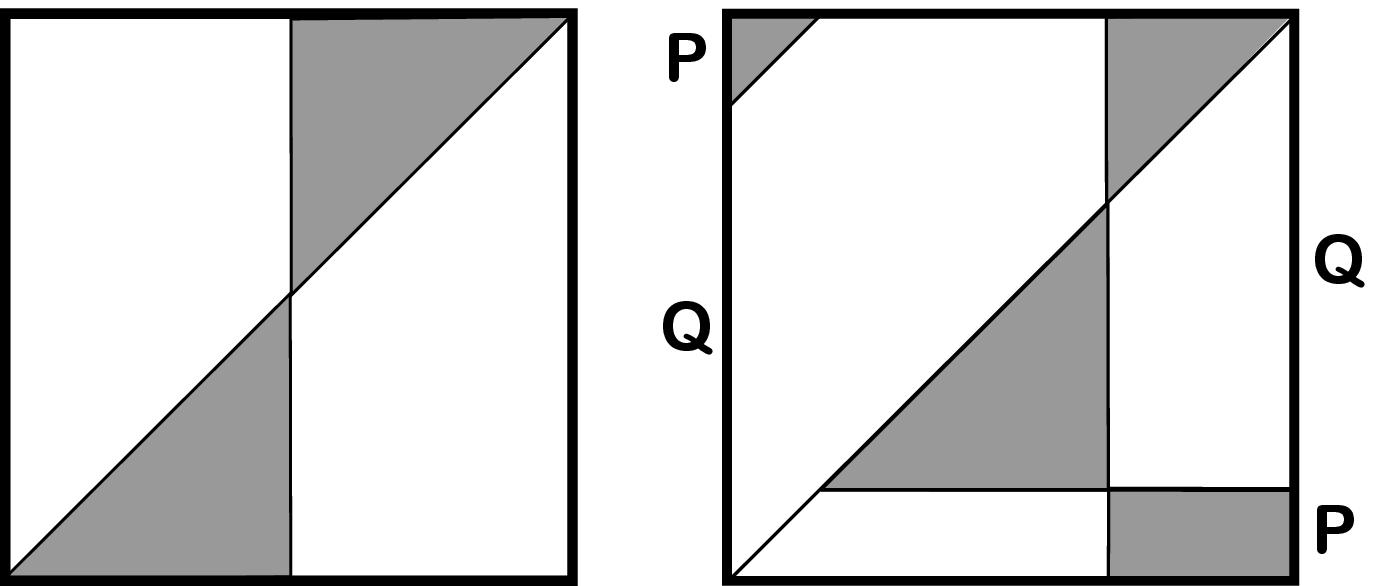}}
\newline
{\bf Figure 6.2:\/} $G'$ on the left and
$G$ on the right.
\end{center}

Combining Lemma \ref{QVERT2} with the description
in this section, we get the following result.

\begin{lemma}
\label{QVERT3}
Let $G=\Psi(G')$.  Then
$\Xi({\cal V\/}_0) \subset \Delta X-G$ and
$\Xi({\cal V\/}_1) \subset G$ and
\end{lemma}

\noindent
{\bf Remark:\/} 
The left and right sides of $G$ really
do match up, because these sides are
identified mod $\Lambda_P$.  In the picture,
the point $(-1,y)$ matches with $(1,y+P)$ if $x<1-P$
and otherwise with $(1,y-Q)$.

\subsection{The Horizontal Case}

We will get a nicer answer if we use the
alternate map
\begin{equation}
\Xi'=\Psi_0^{-1} \circ \Xi, \hskip 30 pt
\Psi_0(T,U_1,U_2)=(T,U_1-P/2-1,U_2-P/2).
\end{equation}
$\Xi'$ differs
from $\Xi$ only by a fiber-preserving
translation.
We define
\begin{equation}
(T',U_1',U_2')=\Psi_0^{-1}(T,U_1,U_2).
\end{equation}
Let $U$ denote the fiber over $T=-1$.

\begin{lemma}
Let $G'$ denote the set in
$U$ such that $(U_1,U_2)$ is a good pair.
Then $\Xi'({\cal H\/}'_0) \subset U-G'$ and
$\Xi'({\cal H\/}'_1) \subset G'$.
\end{lemma}

\startproof
Let $c=(x,y)$ be some point of
$\cal H'_k$.  Using the fact that
$x=2k\omega+\omega/2$, we compute
\begin{equation}
F_P(c)=F_Q(c)=[Py + PQx +1]_2, \hskip 30 pt
F_H(c)=[2Py]_2.
\end{equation}
Let $a_1$ and $a_2$ be these two numbers.
Next, we compute
\begin{equation}
U_1'=[-Py+PQx+1]_2, \hskip 30 pt
U_2'=[Py+PQx]_2.
\end{equation}
Let $b_1$ and $b_2$ be these two numbers.
The numbers $a_1,a_2,b_1,b_2$ satisfy the
hypotheses of Lemma \ref{SWAP}.
Lemma \ref{SWAP} now finishes this proof.
\endproof

\subsection{The Horizontal Case Translated}

Define
\begin{equation}
\rho(x,y)=(x,y+1/2).
\end{equation}
\begin{equation}
\Psi_1(T,U_1.U_2)=
(T,U_1-P/2,U_1+P/2)
\end{equation}

We have
\begin{equation}
\Xi \circ \rho = \rho \circ \Psi_1,
\hskip 30 pt
\rho({\cal H\/}'_k)={\cal H\/}_k.
\end{equation}
Putting all this together, we have
\begin{equation}
\Xi({\cal H\/}_k)=
\Psi_1(\Xi({\cal H\/}'_k))=
\Psi_1\Psi_0\Xi'({\cal H\/}'_k).
\end{equation}
Note that
\begin{equation}
\Psi_1\Psi_0(T,U_1,U_2)=(T,U_1+(1-P),U_2).
\end{equation}
So, $\Xi({\cal H\/}_k)$ is obtained
from $\Xi({\cal H\/}'_k)$ by applying
a horizontal translation of the fiber
by $(1-P,0)$.  Figure 6.3 shows the
sets $G'$ and $G$.

\begin{center}
\resizebox{!}{2.2in}{\includegraphics{Pix/checker6.eps}}
\newline
{\bf Figure 6.3:\/} $G'$ on the left and $G$ on the right.
\end{center}

Define
\begin{equation}
\Psi(U,T_1,T_2)=(U,T_1+1-P,T_2).
\end{equation}

\begin{lemma}
\label{HOR}
Let $G=\Psi(G')$.  Then
$\Xi({\cal H\/}_0) \subset U-G$ and
$\Xi({\cal H\/}_1) \subset G$.
\end{lemma}

\newpage

\section{The Isomorphism Theorem}

\subsection{Nine Sets of Centers}

We want to prove that the grid light set and
the tile light set coincide for any square
tile, and edge of the square, and any parameter.
Using the 
fact that reflections in the horizontal
and vertical midlines of the fundamental 
domain are symmetries of both the
grid light set and the tile light set,
it suffices to prove our result for
the west and south edges of the squares.
We fix some even rational parameter $p/q$.

Recall that ${\cal C\/}$ is the set of centers
of the square tiles.
For each symbol $A \in \{S,W\}$ and
each symbol $B \in \{P,Q\}$ and each
integer $k=0,1,2,...$ we define
${\cal S\/}(A,B,k)$ to be the set 
$c \in {\cal C\/}$ such that the grid description
assigns $k$ light
points of type $B$ to the $A$ edge of
the square centered at $c$.

Here are some easy observations.
\begin{itemize}
\item Since the $\cal P$ lines have slope $-P \in (-1,0)$,
and intersect the $y$-axis at integer points,
${\cal S\/}(H,P,k)=\emptyset$ for $k>1$.

\item Since the $\cal Q$ lines have slope $-Q \in (-2,-1)$
and intersect the $y$-axis at integer points,
${\cal S\/}(H,P,k)=\emptyset$ for $k>2$.

\item Since the $\cal P$ and $\cal Q$ lines intersect
the $y$ axis at integer points, the light
points of  either type on a $\cal V$ line are spaced
integer distances apart.  Hence
${\cal S\/}(V,B,k)=\emptyset$ for $B \in \{P,Q\}$
and $k>1$.
\end{itemize}
So, there are just $9$ nontrivial sets.

As we have mentioned already, our
strategy for understanding the locations
of these sets in the classifying space is
to observe that each one of them contains
one of the images we considered in 
\S \ref{instance}:
\begin{equation}
\label{contain1}
{\cal V\/}_k \subset {\cal S\/}(P,V,k) \cap
{\cal S\/}(Q,V,k), \hskip 30 pt k=0,1.
\end{equation}
\begin{equation}
{\cal H\/}_0 \subset {\cal S\/}(P,V,0) \cap
{\cal S\/}(Q,V,0), \hskip 15 pt
{\cal H\/}_1 \subset {\cal S\/}(P,H,1) \cap
\bigcup_{k=1}^2 {\cal S\/}(Q,H,k).
\end{equation}
Each
of the sets considered in \S \ref{instance}
contains the image of a west or south
edge containing at least one instance of each
particle.  Using the results from
particles obtained in \S \ref{PARTICLE}
we will be able to pin down the images of
the $9$ sets exactly.

\subsection{The West Edges}

We first describe a certain subset
$\Sigma(W,Q,1) \subset X$, fiber by fiber.
\begin{itemize}
\item When the fiber lies over 
$T \in (-1+P,P)$, the set $\Sigma(W,Q,1)$ consists
of all points in the fiber lying on
horizontal lines through the square labeled W.
\item Otherwise, 
$\Sigma(W,Q,1)$ consists
of all points in the fiber lying on
horizontal lines which avoid $W$.
\end{itemize}
Figure 7.1 shows a picture of the two cases.

\begin{center}
\resizebox{!}{1.7in}{\includegraphics{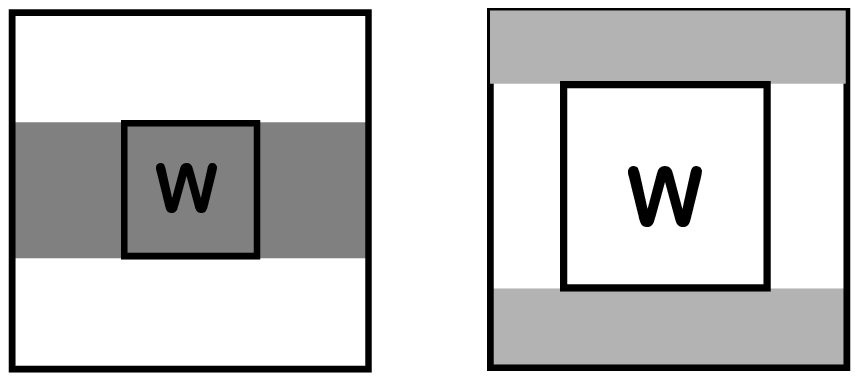}}
\newline
{\bf Figure 7.1:\/} The set
$\Sigma(W,Q,1)$ intersecting a fiber.
\end{center}

The left hand side of Figure 7.2 gives a
further guide to $\Sigma(W,Q,1)$.
We have translated the square torus so that
the W square is at the center of the picture,
but in most fibers the picture will not be
centered this way.

\begin{center}
\resizebox{!}{2.5in}{\includegraphics{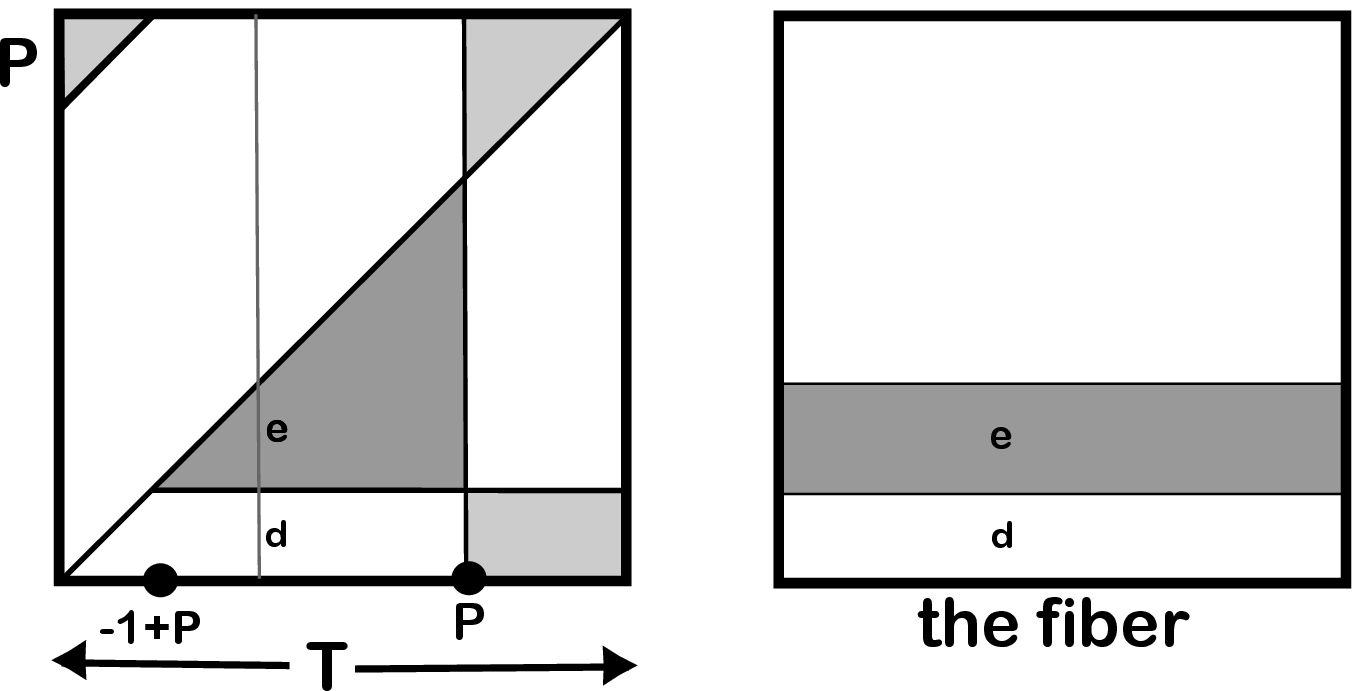}}
\newline
{\bf Figure 7.2:\/} A guide to
understanding $\Sigma(W,Q,1)$.
\end{center}

The picture is drawn when $P \approx 1/3$ but
the set is the same for all choices of $P$,
up to translation (on the torus).
The way the picture works is that you
draw a vertical line $T$ units over from the
left.  You see how this vertical line
intersects the shaded region in the picture.
Then you construct a horizontal strip in the fiber
$[-1,1]$ so that the vertical line intersects
the strip the same way.

The set is not defined in the fiber
$T=P$. However, as discussed in 
\S \ref{defined}, the classifing map
does not map any points of $\cal C$
into this transitional fiber.
Next, define
\begin{itemize}
\item $\Sigma(W,Q,0)=X-\Sigma(W,Q,1)$.
\item $\Sigma(W,P,k)=\Xi(\Sigma(W,P,k)$ for $k=0,1$.
\end{itemize}

\begin{lemma}
\label{tech1}
$\Xi$ maps ${\cal S\/}(V,Q,k)$ into
$\Sigma(V,Q,k)$ for $k=0,1$.
\end{lemma}

\startproof
Comparing Figure 7.2 wish our description of
the set $G \subset \Delta X$ given in
\S \ref{vtrans} and shown on
the right hand side of 
Figure 6.2, we see that
$\Sigma(V,Q,1)$ has the following
description.  We take the union of
the lines parallel to $(0,1,0)$ through
all the points of $G$ and intersect them
with $[-1,1]^3$.  At the same
time we define $\Sigma(V,Q,0)$ has the
same description, using $\Delta X-G$
in place of $G$.
This lemma now follows from
Lemma \ref{QVERT}, 
Lemma \ref{QVERT3} and
Equation \ref{contain1}.
\endproof

\begin{lemma}
\label{tech2}
$\Xi$ maps ${\cal S\/}(V,P,k)$ into
$\Sigma(V,P,k)$ for $k=0,1$.
\end{lemma}

\startproof
This follows from Lemma \ref{tech1},
the definition of $\Sigma(W,P,k)$, and the
symmetry established in Lemma \ref{REFLECT} and
Lemma \ref{refly}.
\endproof

Let $S_1 \Delta S_2$ denote the symmetric
difference of the sets, namely
$$S_1 \Delta S_2=(S_1-S_2) \cup (S_2 - S_2).$$
Looking at the picture fiber by fiber, we see that
\begin{equation}
\label{symm}
\Sigma(W,P,1) \Delta \Sigma(W,Q,1) = X_{WS} \cup X_{WE} \cup X_{WN}.
\end{equation}

Here $X_{WS}$, etc., are the sets which
comprise the partition of the classifying space $X$.
Let $c \in \cal C$ be
a center of a square tile. Let $w$ be the
west edge of the tile centered at $c$.
Suppose that the grid model assigns $0$ or $2$ light
point to $w$.  Then either
$c \in {\cal S\/}(V,B,0)$ for
both choices of $B \in \{P,Q\}$ or
for neither choice.
But then
$$\Xi(c) \not \in \Sigma(V,P,1) \Delta \Sigma(V,Q,1).$$
Hence
$\Xi(c) \not \in X_{WS} \cup X_{WE} \cup X_{WN}$.
But then the tile description assigns a tile to $c$ whose
connector does not involve the $w$.

Suppose on the other hand that $w$ has one light
point.  Then
$$\Xi(c) \in \Sigma(V,P,1) \Delta \Sigma(V,Q,1).$$
Hence
$\Xi(c) \in X_{WS} \cup X_{WE} \cup X_{WN}$,
and the tile description assigns a tile to $c$ whose
connector does involve $w$.
This completes the proof of the
Isomorphism Theorem for the west edges. 
Again, the case of the east edges follows from
symmetry.
\newline
\newline
{\bf Remark:\/}
In writing the equations above, we have been a bit
sloppy about the boundaries of our sets, but
the same analysis as given in 
\S \ref{defined} shows that the image
of $\cal C$ avoids the boundaries of
the relevant sets.

\subsection{The South Edges}

Now we define sets $\Sigma(S,P,k)$ for
$k=0,1$ and $\Sigma(S,Q,k)$ for
$k=0,1,2$. Again, we give a
fiberwise description of these sets.
The picture we show is a diagram
which requires some interpretation.
We will show the picture first
and then give the interpretation.

The left hand side of Figure 7.2 shows
the labeling we use for fibers lying
above points $T \in [1-P,-1+P]$
and the right hand picture shows the
labeing we use for the fibers lying
over the remaining points. Notice
that this is the same dichotomy
that appears in Lemmas \ref{typeP}
and \ref{typeQ}, which deal with
the images of horizontal particles
of type P and Q respectively.

We have drawn the
picture with the square labeled $S$
in the center of the torus, but in
each fiber we mean to construct these
sets based on the true location of $S$.
The diagonal lines are such that they bisect
the edges of the square fundamental domain.

For each $B \in \{P,Q\}$ we define
$\Sigma(S,B,1)$ to be the union, taken
over all fibers, of the regions which have
a $B$ label in them.  The exception is
the label Q2, which does not count as
an instance of the Q label.

\begin{center}
\resizebox{!}{2.8in}{\includegraphics{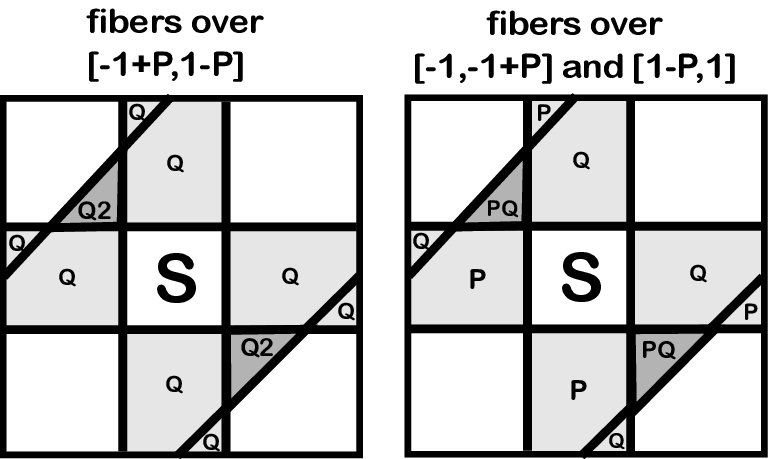}}
\newline
{\bf Figure 7.2:\/} Guide to the sets
\end{center}

Notice also that
$\Sigma(S,P,1) \cap \Sigma(S,Q,1)$ is the
union of all the regions labeled $PQ$.
Finally we define
$\Sigma(S,Q,2)$ as the union
of the regions labeled Q2.
We define $\Sigma(S,P,0)$ each to
be the union of the unlabeled pieces - i.e., the
complement of the sets defined by the labels.

Below we will prove the following result.
\begin{lemma}
\label{TECH}
For each $B \in \{P,Q\}$ and
each relevant $k \in \{0,1,2\}$,
the set ${\cal S\/}(H,B,k)$ is mapped
into $\Sigma(H,B,k)$ by $\Xi$.
\end{lemma}

Assume Lemma \ref{TECH} for the moment.
By inspecting each of the two types of
fibers, we see the following relations
\begin{equation}
\Sigma(S,P,1) \Delta \Sigma(S,Q,1)=X_{SE} \cup X_{SW} \cup X_{SN}.
\end{equation}
\begin{equation}
\Sigma(S,Q,k) \cap (X_{SE} \cup X_{SW} \cup X_{SN})=\emptyset,
\hskip 30 pt k=0,2.
\end{equation}

Assuming Lemma \ref{TECH} and the above
two relations, the rest of the proof of
the Isomorphism Theorem for the south
edges is exactly as it was for the
west edges.  

The rest of the section is devoted to
proving Lemma \ref{TECH}.  We will
break this result up into several
smaller ones.

\begin{lemma}
\label{tech3}
$\Xi$ maps ${\cal S\/}(H,P,k)$ into
$\Sigma(H,P,k)$
for $k=0,1$.
\end{lemma}

\startproof
Looking at Equation \ref{zone1} we see that
as $T \to -1$, the square labeled S shrinks to
the point $(1-P,-1)$.  Therefore, the right
hand side of Figure 6.3, with the shaded
region labeled PQ, shows what our
sets $\Sigma(H,P,k)$ look like above the fiber $T=-1$.
We call these intersections $\Sigma_0(H,P,k)$.
For convenience, we repeat the picture below.

Let $\Sigma_t(H,P,k)$ be the
intersection of $\Sigma(H,P,k)$ with the fiber
over $-1+t$.  We imagine $t \in [-P,P]$ as a
parameter.  We know the picture at $t=0$ and
we want to understand what happens as $t$
varies.

\begin{center}
\resizebox{!}{2in}{\includegraphics{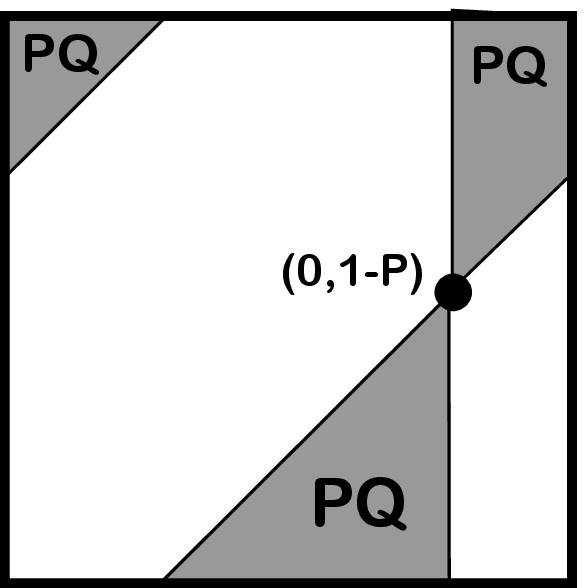}}
\newline
{\bf Figure 7.3:\/} The set $\Sigma_0(H,P,1)$ is shaded.
\end{center}

Looking at the right hand side of Figure 7.2, we see
$\Sigma_t(P,H,k)$ is just a translate of
$\Sigma_0(P,H,k)$.   
Here we have included the portion
of the right hand side of Figure 7.2 which just shows
the regions with the P labels.

\begin{center}
\resizebox{!}{2.2in}{\includegraphics{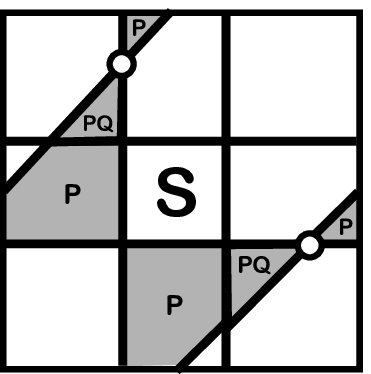}}
\newline
{\bf Figure 7.4:\/} $\Sigma(H,P,1)$ over a fiber
in $[-1+P,1] \cup [1-P,1]$.
\end{center}

Call the marked points in Figure 7.4 the
{\it apices\/} of the fiber.  Each apex lies on
a vertical geodesic extending a
vertical edge of the square $S$ and is half way
around the fiber from the center of this edge.
From this description, and from Equations
\ref{zone1} and \ref{zone3} (and, to be
honest, from computer drawings) we see that
the apices move with speed $1$ along lines
of slope $1$ as a function of $t$.  In
otherwise, we obtain $\Sigma_t(H,P,k)$ from
$\Sigma_0(H,P,k)$ by translating
along the vector $(t,t,t)$.

Combining this information with Lemma \ref{typeP} and
Lemma \ref{HOR}, we get the conclusion of this lemma.
\endproof

\begin{lemma}
\label{tech4}
Let $c \in {\cal S\/}(H,Q,k)$ for
some $k$ and suppose that
$\Xi(c)$ is contained in a fiber
over $T \in [-1,-1+P]$.  Then
$\Xi(c) \in {\cal S\/}(H,Q,k)$.
\end{lemma}

\startproof
The proof is similar to what we did in
Lemma \ref{tech4}. 
For $t \in [-1,-1+P]$ the
set $\Sigma_t(H,Q,2)$ is empty, and
the set $\Sigma_t(H,Q,1)$ is a translate
of the set $G$ from Lemma \ref{HOR}. 
Each apex lies on
a horizontal geodesic extending a
horizontal edge of the square $S$ and is half way
around the fiber from the center of this edge.
From this description, and from Equations
\ref{zone1}, we see that the apices do not
move at all as $t$ varies.
In other words, we obtain
$\Sigma_t(H,Q,k)$ from
from $\Sigma_{-1}(H,Q,k)$ by applying 
the translation $(t+1,0,0)$.
It now follows from Lemma
\ref{typeQ} and Lemma \ref{HOR} that
$\Xi$ maps each $c \in {\cal S\/}(H,Q,k)$ into
$\Sigma(H,Q,k)$ for $k=1,2$ provided that
$\Xi(c)$ lies in a fiber over $[-1+P,P]$.
It follows from Lemma \ref{typeQ} that
$\Xi({\cal S\/}(H,Q,2)$ does not intersect
these fibers at all.
\endproof

\begin{lemma}
\label{tech5}
Let $c \in {\cal S\/}(H,Q,k)$ for
some $k$ and suppose that
$\Xi(c)$ is contained in a fiber
over $T \in [1-P,1]$.  Then
$\Xi(c) \in {\cal S\/}(H,Q,k)$.
\end{lemma}

\startproof
For $t \in [1-P,1]$ the analysis is exactly
the same, except that $\Sigma_t(H,Q,k)$ is
translated in the fiber by $(P,P)$.  
Were we to have chosen the branch of $[\cdot]_2$
that took values in $(-2,2]$ rather than 
in $[-2,2)$ we would have proved a version
of Lemma \ref{HOR} for sets in the fiber
over $T=1$, and our set $G$ would have
been translated by $(P,P)$, thanks to the
boundary identifications on the fundamental
domain boundary coming from the lattice
$\Lambda_P$. 
\endproof

\begin{lemma}
\label{tech6}
Let $c \in {\cal S\/}(H,Q,k)$ for
some $k$ and suppose that
$\Xi(c)$ is contained in a fiber
over $T \in [-1+P,1-P]$.  Then
$\Xi(c) \in {\cal S\/}(H,Q,k)$.
\end{lemma}

\startproof
Let $\pi: X \to \R^2$ denote the
projection $\pi(T,U_1,U_2)=(U_1,U_2)$.
For $t \in [-1+P,1-P]$ we have

\begin{equation}
\label{final0}
\pi(\Sigma_t(H,Q,0))=\pi(\Sigma_{-1}(H,Q,0) \cap \pi(\Sigma_1(H,Q,0)).
\end{equation}

\begin{equation}
\label{final1}
\pi(\Sigma_t(H,Q,1))=\pi(\Sigma_{-1}(H,Q,1) \Delta \pi(\Sigma_1(H,Q,1)).
\end{equation}

\begin{equation}
\label{final2}
\pi(\Sigma_t(H,Q,2))=\pi(\Sigma_{-1}(H,Q,1) \cap \pi(\Sigma_1(H,Q,1)).
\end{equation}

Here $\Delta$ denotes symmetric difference.  Figure 7.4 illustrates
these equalities.  We have taken the left side of
Figure 7.2 and re-coloring it so that
the different shades show the different
translates of the basic set $G$ from
Lemma \ref{HOR}.
The translation $(U_1,U_2) \to (U_1+P,U_2+P)$
carries the white point in Figure 7.4 to the
black point.  In Figure 7.4,
$\Sigma_t(H,Q,1)$ is shown in two different
shades of grey and
$\Sigma_t(H,Q,2)$ in black.
The white set is 
$\Sigma_t(H,Q,0)$.

\begin{center}
\resizebox{!}{3in}{\includegraphics{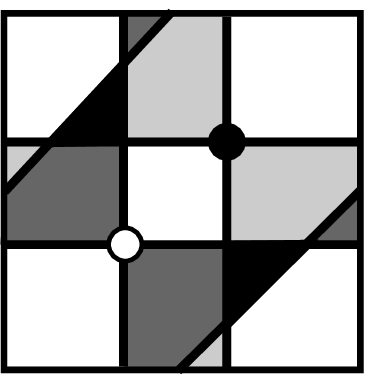}}
\newline
{\bf Figure 7.4:\/} $\Sigma(H,P,1)$ (grey) and
$\Sigma(H,P,2)$ over a fiber in
$[-1+P,1-P]$.
\end{center}

Now there are $4$ cases to consider.
\newline
\newline
{\bf Case 1:\/}
Suppose that the south edge of $c$ contains two
light points, $z_1$ and $z_2$.
Suppose that $z_1$ is the one on the left and
$z_2$ is the one on the right.
By Lemma \ref{typeQQ},
there is another instance $z_2'$ of the particle
containing $z_2$ such that
$\Xi(z_2') \in \Sigma_{-1}(H,Q,1)$, and
$\Xi(z_2')$ and $\Xi(z_2)$ lie on the same
horizontal segment in $[-1,1]^3$. 
Similarly, there is another instance $z_1'$
of the particle containing $z_1$ so that
$\Xi(z_1') \in \Sigma_1(H,Q,1)$, and
$\Xi(z_1')$ and $\Xi(z_1)$ lie on the same
horizontal segment in $[-1,1]^3$.
But then it follows from Equation
\ref{final2} that
$\Xi(c) \in \Sigma(H,Q,2)$.
\newline
\newline
{\bf Case 2:\/}
Suppose that the south edge of $c$ contains one
light point $z_1$ and one dark point $z_2$, with
$z_1$ being on the left.  This is the
same as Case 1, except that
$\Xi(z_2') \in \Sigma_{-1}(H,Q,0)$.
But then, by Equation \ref{final1}, we
have $\Xi(c) \in \Sigma(H,Q,1)$.
\newline
\newline
{\bf Case 3:\/}
Suppose that the south edge of $c$ contains one
dark point $z_1$ and one light point $z_2$, with
$z_1$ being on the left. This is the
same as Case 1, except that
$\Xi(z_1') \in \Sigma_{1}(H,Q,0)$.
But then, by Equation \ref{final1}, we
have $\Xi(c) \in \Sigma(H,Q,1)$.
\newline
\newline
{\bf Case 4:\/}
Suppose that the south edge of $c$ contains 
two dark points.  This time we have
$\Xi(z_1') \in \Sigma_{1}(H,Q,0)$, and
$\Xi(z_2') \in \Sigma_{-1}(H,Q,0)$.
Equation \ref{final0} tells us that
$\Xi(c) \in \Sigma(H,Q,0)$.
\newline

This exhausts all the cases.
\endproof

This concludes the proof of the Isomorphism Theorem.

\newpage

\section{The PET Equivalence Theorem}

\subsection{Orientations}
\label{orient}

In this chapter we deduce the PET Equivalence
Theorem from the Isomorphism Theorem.
This is mainly just a matter of
interpretation.

The tile description of the plaid model
allows for the possibility of orienting
the loops in the model.  There are $12+1$
ways of drawing a directed edge in a square
tile which connects up the midpoints of
different sides.  Aside from the empty
tile, the remaining $12$ tiles are indexed
by ordered $2$-element subsets of
$\{N,S,E,W\}$.   As a first step
in proving the PET Equivalence Theorem.
we explain the oriented version of
our model.

We partition $[0,1] \times \R^3$ by
thinking of this space as the universal
cover of $X$ and then lifting the pieces of
the partition.  We then quotient out by
the index $2$ affine subgroup $\Lambda'$ which,
for each parameter $P$, acts as the lattice
generated by
\begin{equation}
(4,2P,2P), \hskip 30 pt
(0,2,0), \hskip 30 pt (0,0,2).
\end{equation}
We let $\widehat X$ denote this quotient.
We think of 
$$[0,1] \times [-2,2] \times [-1,1]^2$$
as the fundamental domain for $\widehat X$.
This turns out to be the most symmetric
choice.  We can identify $X$ with the
middle half of this fundamental domain,
namely
$$[0,1] \times [-1,1] \times [-1,1]^2,$$

We enhance the labels on $\widehat X$
in the following way:
\begin{itemize}
\item For $T \in [-1,1]$ we use the labels as they
are ordered in Figure 3.1.
\item For $T \in [-2,-1] \cup [1,2]$ we use the labels 
in the reverse order.  So, for instance, the piece
labeled NW in Figure 3.1 would be labeled WN in
the other half of $\widehat X$.
\end{itemize}

One can see the oriented version in action using
my computer program.

\subsection{Proof of Oriented Coherence}

\begin{theorem}
The classifying pair
$(\widehat \Xi_P,\widehat X_P)$ produces
coherently oriented tilings at all
parameters $P \in (0,1)$.
\end{theorem}

We will prove this result through a
series of smaller lemmas.
We will explain the
argument in one of $8$ cases.  The argument works
the same in all cases.  Define
\begin{equation}
\widehat X_{S \downarrow}=\widehat X_{ES} \cup \widehat X_{WS} \cup
\widehat X_{NS}
\end{equation}
as the union of all regions in $\widehat X$ which assign
tiles that point into their south edge.  This set
is a finite union of $4$ dimensional convex integral polytopes.
In a similar way, define
$\widehat X_{N,\downarrow}$ as the set of all regions which
assign tiles that point out of their north edges.

There is an integral affine transformation
$\Upsilon: \widehat X \to \widehat X$ which has the property that
\begin{equation}
\Upsilon \circ \Xi(x,y)=\Xi(x,y-1).
\end{equation}
From Equation \ref{map}, we get
\begin{equation}
\Upsilon(P,\widehat T,U_1,U_2)=(P,\widehat T-2,U_1,U_2-2P) \hskip 10 pt
{\rm mod\/} \hskip 10 pt \Lambda.
\end{equation}
This map is a fiberwise translation.
In particular, the image of a finite union of $K$ convex integer
polytopes under $\Upsilon$ is again a finite union of $K$
convex integer polytopes.  

The statement that every tile which points into the south edge
meets a tile below it that points out of the north edge is
equivalent to the statement that
\begin{equation}
\label{mesh}
\Upsilon(\widehat X_{S,\downarrow})=\widehat X_{N,\downarrow}.
\end{equation}

Recall from the description in \S \ref{partition0}
that $X$ is partitioned into $3$ zones.
So, $\widehat X$ is partitioned into $6$ zones.
More precisely, $\widehat X$ is a fiber bundle
over a union of $6$ integer triangles, and
each zone is the preimage of one of the triangles.
Each of our sets intersects a zone in a
convex integer polytope.  Hence $K=6$.

Say that a {\it triad\/} is a triple of
parameters $\{(P_i,T_i)\}$ for $i=1,2,3$ which
lies in the interior of one of the $6$ triangles
in the base space.

\begin{lemma}
Suppose that Equation \ref{mesh} holds in
the fibers associated to $6$ triads,
one per zone.  Then Equation \ref{mesh}
holds everywhere.
\end{lemma}

\startproof
The two sets in Equation \ref{mesh} intersect
each zone in a convex integer polytope.
Let $L$ and $R$ respectively be the intersection
of the left and right hand sides of
Equation \ref{mesh} with one of the zones.
As we discussed in \S \ref{partition0}, both
$L$ and $R$ are intersect
each interior fiber in a rectangle.
(We mean a fiber over an interior point of the base.)
From this description, 
we have $L=R$ provided that $L$ and $R$ intersect
a triad of fibers in the same way.
\endproof

\begin{lemma}
Suppose that the classifying pair
$(\widehat \Xi_P,\widehat X)$ produces
a coherently oriented tiling
at the parameters $p/q=3/8$ and
$p=4/11$ (corresponding to $P=6/11$ and
$P=8/15$).  Then Equation \ref{mesh}
holds everywhere.
\end{lemma}

\startproof
First we consider some general parameter $p/q$.
Let $\omega=p+q$.  Let $P=2p/\omega$.
Let $\cal C$ denote the set of centers of
the square tiles.
Over any fiber $F$ which contains points
of $\Xi_P({\cal C\/})$, we divide
$F$ into $\omega^2$ squares of side-length
$2/\omega$.  As we discussed in
\S \ref{defined}, the image of
$\Xi_P({\cal C\/}) \cap F$ is exactly
the set of centers of these little squares.

At the same time, 
$\widehat X_{N,\downarrow}$ intersects $F$
in the segments extending the edges of the
little squares.    But the translation
$\Upsilon$ preserves the division of
$\widehat X$ into cubes of side length
$2/\omega$.  Hence $\Upsilon(\widehat X_{S,\downarrow})$
also intersects $F$ in the segments extending the
sides of the little squares.

It follows from Lemma \ref{bijection} that
$\Xi({\cal C\/}) \cap F$ contains every
little square center.  Hence, of our two
sets disagreed in $F$, there would be some
mismatch of orientations we could see in
the tiling.  Since the tiling is coherently
oriented, there is no mismatch.

The parameters $3/8$ and $4/11$ are large
enough so that the image
$$\Xi_{6/11}({\cal C\/}) \cup
\Xi_{8/15}({\cal C\/})$$
contains a triad in every zone.
\endproof

Now we simply inspect the picture for these
two parameters and see that the tiling
produced is coherent.  The parameters involved
are so small that there is no round-off
error in the calculation.

\subsection{Polytope Exchange Transformations}

The data for a polytope exchange transformation is
two partitions of a polytope $Y$ into smaller
polytopes.
\begin{equation}
Y=\bigcup_{i=1}^N L_i=\bigcup_{i=1}^N R_i.
\end{equation}
By {\it partition\/} we mean that the 
smaller pieces have pairwise disjoint interiors.

We insist that there is some translation
$f_i$ so that $f_i(L_i)=R_i$ for all $i$.
These maps piece together to give an almost-everywhere
defined map $F: Y \to Y$.  One defines
$F(p)=f_i(p)$ provided that $p$ lies in the interior
of $L_i$.  The inverse map is defined to be
$F^{-1}(q)=f_i^{-1}(q)$ provided that $q$ lies
in the interior of $R_i$.

A variant of this is to consider $Y$ to
be a flat manifold.  This case really isn't
that different from the polytope case, and indeed
by refining the partitions we can
readily convert between one to the other.
One encounters this situation when one considers
a circle rotation to be a $2$-interval IET.
A refinement of the partition introduces new
points where the map is not, strictly speaking,
defined. However, we will see that this
situation has no impact on our main goals.
In the irrational case, we will study orbits
which are well-defined in either set-up.

For each parameter $p/q$, we set $P=2p/(p+q)$, 
and our space 
is $\widehat X_P$.
The partition is given by
\begin{equation}
\widehat X_{S \downarrow} \cup
\widehat X_{W \leftarrow} \cup
\widehat X_{N \uparrow} \cup
\widehat X_{E \rightarrow}.
\cup \widehat X_{\Box}.
\end{equation}
The remaining set $\widehat X_{\Box}$ is just the complement.
These sets have an obvious meaning.  We already defined
$\widehat X_{S \downarrow}$.  The set
$\widehat X_{E,\rightarrow}$ is the union of
regions which assign tiles which point into their east edges.
And so on.
The second partition of $\widehat X$ is obtained by
reversing all the arrows.

In other words, the pieces in the first partition
are understood as the ones which specify tiles according
to the edges they point into and the pieces in the
second partition are understood as the ones which specify
tiles according to the edges they point out of.  Both
partitions share the regions which assign the empty tile.

We have $5$ {\it curve following maps\/}
$\Upsilon_{\Box}: \widehat X_{\Box} \to \widehat X_{\Box}$ is
the identity, and then
$\Upsilon_S: \widehat X_{S \downarrow} \to \widehat X_{N\downarrow}$.
And so on.
These maps are all translations on their domain.

There is a natural map on $\cal C$, the set of centers
of the square tiles.  We simply follow the directed edge
of the tile centered at $c$
and arrive at the next tile center.  By construction,
 $\widehat \Xi_P$ conjugates this map on
$\cal C$ to our polyhedron exchange transformation.
If we work within the fundamental domain, then
$\widehat \Xi_P$ is a bijection, by Lemma \ref{bijection}.
By construction $\widehat \Xi_P$ sets up a tautological
bijection between the plaid polygons and the
set of orbits of points of the form
$\widehat \Xi(c)$, where $c \in \cal C$.  We call
these orbits {\it special orbits\/}.

By construction $\widehat \Xi_P$
sets up a dynamics-respecting bijection 
the special orbits and the plaid polygons
contained in the fundamental domain, with
respect to the parameter $p/q$.

To get the vector dynamics mentioned in the
PET Equivalence Theorem, we assign the vectors
$(0,-1)$ to the region $\widehat X_{S\downarrow}$, 
and $(-1,0)$ to $\widehat X_{W\leftarrow}$, etc.
When we follow the orbit of
$\Xi(c)$, we get record the list of vectors labeling
the regions successively visited by the point. This gives
vectors $v_1,v_2,v_3$.  The vectors
$c$, $c+v_1$, $c+v_1+v_2$, etc. are the vertices
of the plaid polygon containing $c$.
This is what we mean by saying that the plaid polygons
describe the vector dynamics of the special orbits
of the PET.  This proves the first statement
of the PET Equivalence Theorem.

\subsection{Global Picture}

Now we turn to the second statement of the
PET Equivalence Theorem.
We define an {\it affine polytope exchange\/}
transformation just as above, except that the
maps $f_i$ are affine transformations
rather than translations.  We say that such
dynamical system is {\it fibered\/} if it
is a fiber bundle over some base space such
that the maps preserve each fiber and act
as ordinary PETs there.  

Our
total space $\widehat X$ is naturally a flat affine
manifold, as we discussed (for the two-fold
quotient) in \S 3.
The partitions defined in the previous section piece
together to give a partition of $\widehat X$ into
finitely many convex polytopes.  The polyhedron
exchange maps piece together to give piecewise
affine maps of $\widehat X$.  The definitions
make sense even for irrational parameters.
Thus $\widehat X$ is a fibered piecewise affine
polytope exchange transformation, and the fiber
over $P$ is $\widehat X_P$.

We have
\begin{equation}
\widehat X=([0,1] \times \R^3)/\widehat \Lambda,
\end{equation}
where $\widehat \Lambda$ is a group of
{\it integral\/} affine transformations.
We say that a convex polytope in $\widehat X$ is
{\it integral\/} if one of its lifts to
the universal cover has integer vertices.
This definition is independent of lift.
By construction, all the polytopes in our
partitions are integral.
In short, $\widehat X$ is a fibered, integral,
piecewise affine polytope exchange transformation.
This proves the second statement of the
PET Equivalence Theorem.

\newpage

\section{Connection to Outer Billiards}
\label{ob}

\subsection{Basic Definitions}

Outer billiards is a dynamical system defined
on the outside of a compact convex set $K$.  The map is
defined like this.  One starts with $p_0 \in \R^2-P$
and then defines $p_1$ so that the line segment
$\overline{p_0p_1}$ is tangent to $K$ at a vertex with
bisects this segment, as shown in Figure 9.1.

\begin{center}
\resizebox{!}{2.2in}{\includegraphics{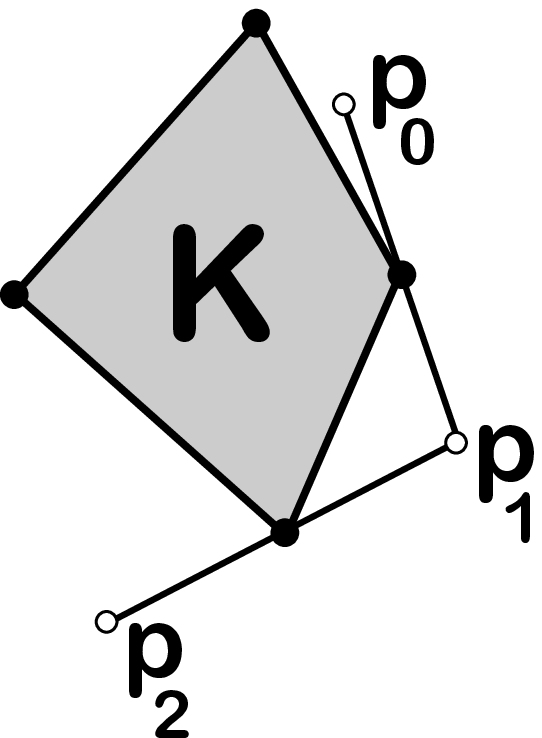}}
\newline
{\bf Figure 9.1:\/} Outer billiards on a kite
\end{center}

Of two possible choices for $p_1$, we choose so that the
one finds $K$ on the right as one walks from $p_0$ to $p_1$.
One defines $p_2,p_3,...$ in the same way.  The map
is not defined when $p_0$ lies on a finite union of
rays. This implies that the full orbit $\{p_n\}_{n \in \Z}$
is defined for all $p_0$ in the complement of a countable
union of line segments.  

Outer billiards was introduced by B. H. Neumann in the
late $1950$s. See [{\bf N\/}].
One of the central questions about outer billiards has
been: Does there exist
a convex set $K$ (not necessarily a polygon) and a point
$p_0$ so that the orbit $\{p_n\}$ exits every compact
set?  This question, dating from around $1960$, is called
the Moser-Neumann question.  See  [{\bf M\/}].
  My monograph [{\bf S1\/}] has a
survey of the known results about this problem. 
See also [{\bf T1\/}] and [{\bf T2\/}]. The main
work on this Moser-Neumann problem is contained in 
[{\bf D\/}], [{\bf VS\/}], [{\bf K\/}], [{\bf GS\/}], [{\bf G\/}],
[{\bf S1\/}], [{\bf S2\/}], and [{\bf DF\/}].

In [{\bf S2\/}], I found the first example of a shape
with unbounded orbits - the Penrose kite.
In [{\bf S1\/}] I subsequently showed that
there are orbits with respect to any irrational kite.
Any kite is affinely equivalent to the kite
$K_A$ with vertices
\begin{equation}
(-1,0), \hskip 30 pt
(0,1), \hskip 30 pt (0,-1), \hskip 30 pt (A,0),
\hskip 30 pt A \in (0,1).
\end{equation}
$K_A$ is rational if and only if $A \in \Q$.
Outer billiards on
$K_A$ preserves the infinite family of
horizontal lines
\begin{equation}
\Lambda=\R \times \Z_1.
\end{equation}
Here $\Z_1$ is the set of odd integers.
Let $\psi_A$ denote the square of the outer billiards
map, restricted to $\Lambda$.

\subsection{The Arithmetic Graph}

One of the main tools I used in [{\bf S1\/}]
to understand $\psi_A$ was the
{\it arithmetic graph\/}.
Let $\Psi_A$ denote
the second iterate of the first return map
of $\psi_A$ to the two lines
\begin{equation}
\Upsilon=\R \times \{-1,1\} \subset \Lambda.
\end{equation}
(The first return map moves points a long
distance, but the second return map is
close to the identity.)
The set $\Upsilon$ is shown in Figure 9.2.
If we choose some ``offset'' $t \in \R_+$, as
well as integers $m_0,n_0$, then it turns out
that
\begin{equation}
\Psi_A\Big(2m_0A+2n_0+t,(-1)^{m_0+n_0}\Big)=
\Big(2m_1A+2n_1+t,(-1)^{m_1+n_1}\Big)
\end{equation}
for some other pair of integers $m_1,n_1$ which
depends on $A$ and $t$.

\begin{center}
\resizebox{!}{2.2in}{\includegraphics{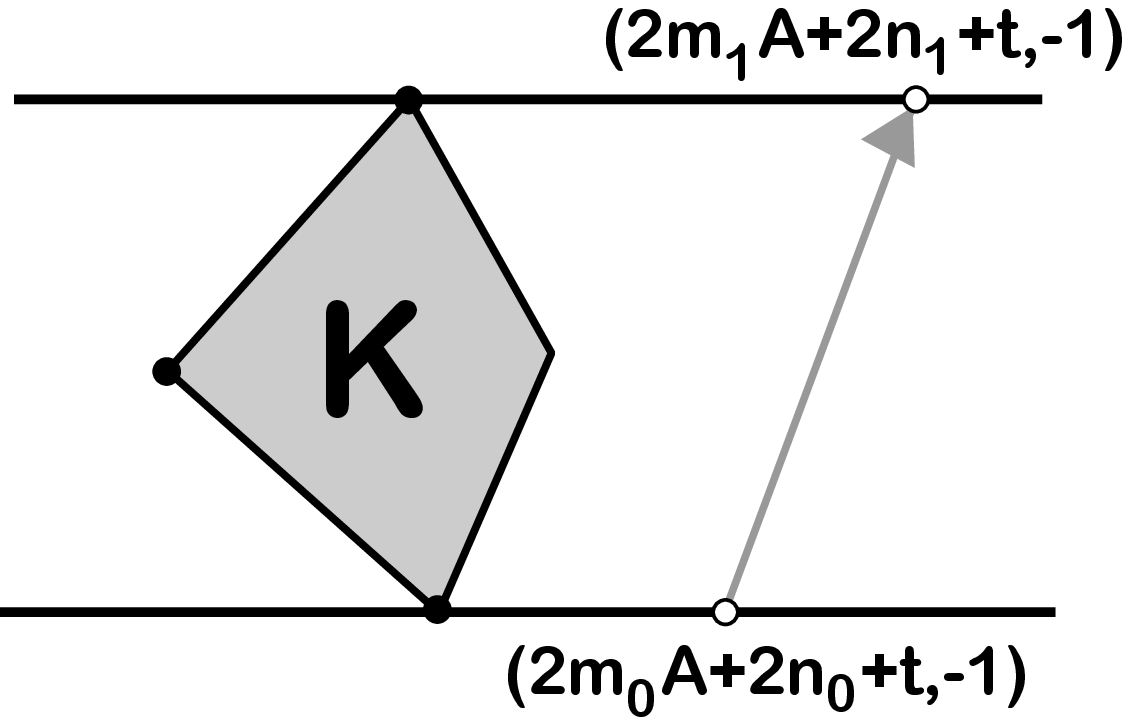}}
\newline
{\bf Figure 9.2:\/} The second return map to $\Upsilon$.
\end{center}

To the point 
\begin{equation}
\label{arithmetic}
p_0=\Big(2m_0A+2n_0+t,(-1)^{m_0+n_0}\Big).
\end{equation}
we have the orbit $\{p_k\}$ and, by the above
correspondence,
the sequence of integers
$\{(m_k,n_k)\}$.
These integers define a polygonal path whose
vertices lie in the integer lattice - i.e. a
{\it lattice path\/}. 
This lattice path is called the
{\it arithmetic graph\/} of the orbit.

I proved in [{\bf S1\/}] that these lattice
paths are always embedded.  It we fix $t$ and
consider the orbits of all possible points
of the form given in Equation \ref{arithmetic},
we get a disjoint union of embedded lattice paths.
When $A=p/q$ is rational, it turns out that
$\Upsilon$ is partitioned into intervals of size
$2/q$ which are permuted by $\Psi_A$.  In this
case, the choice $t=1/q$ is a canonical choice;
it corresponds to the centers of these intervals.
The corresponding union of lattice paths describes
every special orbit at the same time.  We call
this union of all arithmetic graphs the
{\it arithmetic landscape\/}, and we denote it
by $\Gamma(p/q)$. 

\subsection{Quasi-Isomorphism Conjecture}
\label{dyn}

It turns out that $\Gamma(p/q)$ is invariant
under a certain lattice acting on $\Z^2$.  The
lattice is generated by vectors $V$ and
$(p+q)W$, where
\begin{equation}
V=(q,-p), \hskip 30 pt
W=\bigg(\frac{2pq}{p+q},\frac{2pq+q^2-p^2}{p+q}\bigg).
\end{equation}
A fundamental domain for the action is a
certain parallelogram with
sides parallel to $V$ and $(p+q)W$ and southwest
corner $\Theta(p,q)$.  The point $\Theta(p,q)$ is
given by a rather complicated formula, which
runs like this:
\begin{itemize}
\item Let $p'/q'$ be the unique even
rational with $|pq'-p'q|=1$.
\item Let $x=(q-q')/2$ and $y=-xp/q$.
\item If $p'/q'>p/q$ let $x=-x$ and $y=-y$.
\item If $x>0$ let $x=x-q$ and $y=y+0$.
\end{itemize}
Then $\Theta=(x,y)$.

The {\it first block\/} of the big parallelogram
has southwest corner $\Theta$ and sides
parallel to $V$ and $W$.  Call this parallelogram
$\Omega_0$.  It turns out that none of the
polygons in the arithmetic landscape
cross the boundaries of $\Omega_0$.  This
is about half the content of our Hexagrid Theorem
in the case when $pq$ is even. (In fact we
only proved the odd case.)

There is an affine transformation $T_{p/q}$
which maps the first block $B_0=[0,p+q]^2$ to
$\Omega_0$.  Figure 1.3 shows
$T_{p/q}^{-1}(\Gamma_{3/8})$ superimposed
over the union of plaid polygons for the
parameter $3/8$ contained in $B_0$.
One can see similar pictures for other
parameters using my computer program.
We can define the other $p+q-1$ blocks by moving
$\Omega_0$ parallel to itself.  Our affine map
$T_{p/q}$ carries the fundamental domain for the
plaid model to the fundamental domain for the
arithmetic graph.  

\begin{conjecture}[Quasi-Isomorphism]
\label{qi}
There is  a bijection 
between the polygons of $T_{p/q}^{-1}(\Gamma_{p/q})$
in the fundamental domain and the
polygons of the plaid model at $p/q$ in the fundamental
domain, so that each polygon in the one model is within
a $2$-tubular neighborhood of the corresponding
polygon in the other.
\end{conjecture}

One might wonder what the Quasi-Isomorphism
Conjecture has to do with the dynamics.
Here is a consequence of the
Quasi-Isomorphism Conjecture (which we
will flesh out when we prove the conjecture).
Let $A=p/q$ as above.
Recall that $\Psi_A$ is the second return
map to $\Upsilon$.
We define two orbits of $\Psi_A$ to be
{\it equivalent\/} if reflection in the
$x$-axis swaps them.  Two equivalent
orbits are essentially the same.  All the
$\Psi_A$ orbits come in pairs like this.

\begin{corollary}
Let $A=p/q$ with $pq$ even.
Let $\pi_1$ be projection onto the $x$-axis.
There is a bijection between the plaid polygons
in $[0,\infty] \times [0,p+q]$ and the 
equivalence classes of $\Psi_A$ orbits.  
The bijection is such that
for each polygon $\gamma$, the image
$\pi_1(\gamma)$, when parametrized by arc length,
remains vertex by vertex within $5$ units of
the corresponding $\Psi_A$ orbit.
\end{corollary}

In other words, the plaid model says all there is
to know about the coarse geometry and arithmetic of the
special outer billiards orbits on kites.

\newpage

\section{References}

[{\bf DeB\/}] N. E. J. De Bruijn, {\it Algebraic theory of Penrose's nonperiodic tilings\/},
Nederl. Akad. Wentensch. Proc. {\bf 84\/}:39--66 (1981).
\newline
[{\bf DF\/}] D. Dolyopyat and B. Fayad, {\it Unbounded orbits for semicircular
outer billiards\/}, Annales Henri Poincar\'{e} {\bf 10\/} (2009) pp 357-375
\newline
[{\bf G\/}] D. Genin, {\it Regular and Chaotic Dynamics of
Outer Billiards\/}, Pennsylvania State University Ph.D. thesis, State College (2005).
\newline
[{\bf GS\/}] E. Gutkin and N. Simanyi, {\it Dual polygonal
billiard and necklace dynamics\/}, Comm. Math. Phys.
{\bf 143\/}:431--450 (1991).
\newline
[{\bf H\/}] W. Hooper, {\it Renormalization of Polygon Exchange Transformations
arising from Corner Percolation\/}, Invent. Math. {\bf 191.2\/} (2013) pp 255-320
\newline
[{\bf Ko\/}] Kolodziej, {\it The antibilliard outside a polygon\/},
Bull. Pol. Acad Sci. Math.
{\bf 37\/}:163--168 (1994).
\newline
[{\bf M\/}] J. Moser, {\it Is the solar system stable?\/},
Math. Intelligencer {\bf 1\/}:65--71 (1978).
\newline
[{\bf N\/}] B. H. Neumann, {\it Sharing ham and eggs\/},
Summary of a Manchester Mathematics Colloquium, 25 Jan 1959,
published in Iota, the Manchester University Mathematics Students' Journal.
\newline
[{\bf S1\/}] R. E. Schwartz, {\it Outer Billiard on Kites\/},
Annals of Math Studies {\bf 171\/} (2009)
\newline
[{\bf S2\/}] R. E. Schwartz, {\it Unbounded Orbits for Outer Billiards\/},
J. Mod. Dyn. {\bf 3\/}:371--424 (2007). 
\newline
[{\bf T1\/}] S. Tabachnikov, {\it Geometry and billiards\/},
Student Mathematical Library 30,
Amer. Math. Soc. (2005).
\newline
[{\bf T2\/}] S. Tabachnikov, {\it Billiards\/}, Soci\'{e}t\'{e} Math\'{e}matique de France, 
``Panoramas et Syntheses'' 1, 1995
\newline
[{\bf VS\/}] F. Vivaldi and A. Shaidenko, {\it Global stability of a class of discontinuous
dual billiards\/}, Comm. Math. Phys. {\bf 110\/}:625--640 (1987).

\end{document}